
\def\LaTeX{%
  \let\Begin\begin
  \let\End\end
  \let\salta\relax
  \let\finqui\relax
  \let\futuro\relax}

\def\UK{\def\our{our}\let\sz s}
\def\USA{\def\our{or}\let\sz z}

\UK



\LaTeX

\USA


\salta

\documentclass[twoside,12pt]{article}
\setlength{\textheight}{24cm}
\setlength{\textwidth}{16cm}
\setlength{\oddsidemargin}{2mm}
\setlength{\evensidemargin}{2mm}
\setlength{\topmargin}{-15mm}
\parskip2mm


\usepackage[usenames,dvipsnames]{color}
\usepackage{amsmath}
\usepackage{amsthm}
\usepackage{amssymb}
\usepackage[mathcal]{euscript}
\usepackage{enumitem}
\usepackage{cite}
\usepackage{hyperref}

\definecolor{rosso}{rgb}{0.8,0,0}

\def\an #1{#1}



\bibliographystyle{plain}


%

\finqui

\def\Beq{\Begin{equation}}
\def\Eeq{\End{equation}}

\def\Bthm{\Begin{theorem}}
\def\Ethm{\End{theorem}}

\def\Bcor{\Begin{corollary}}
\def\Ecor{\End{corollary}}
\def\Brem{\Begin{remark}\rm}
\def\Erem{\End{remark}}

\def\Bdim{\Begin{proof}}
\def\Edim{\End{proof}}
\def\Bcenter{\Begin{center}}
\def\Ecenter{\End{center}}
\let\non\nonumber




\def\step #1 \par{\medskip\noindent{\bf #1.}\quad}


\def\Lip{Lip\-schitz}
\def\Holder{H\"older}
\def\Frechet{Fr\'echet}

\def\wk{well-known}
\def\socal{so-called}
\def\lhs{left-hand side}
\def\rhs{right-hand side}
\def\sfw{straightforward}


\def\multibold #1{\def\arg{#1}%
  \ifx\arg\pto \let\next\relax
  \else
  \def\next{\expandafter
    \def\csname #1#1#1\endcsname{{\bf #1}}%
    \multibold}%
  \fi \next}

\def\pto{.}

\def\multical #1{\def\arg{#1}%
  \ifx\arg\pto \let\next\relax
  \else
  \def\next{\expandafter
    \def\csname cal#1\endcsname{{\cal #1}}%
    \multical}%
  \fi \next}


\def\multimathop #1 {\def\arg{#1}%
  \ifx\arg\pto \let\next\relax
  \else
  \def\next{\expandafter
    \def\csname #1\endcsname{\mathop{\rm #1}\nolimits}%
    \multimathop}%
  \fi \next}

\multibold
qwertyuiopasdfghjklzxcvbnmQWERTYUIOPASDFGHJKLZXCVBNM.

\multical
QWERTYUIOPASDFGHJKLZXCVBNM.

\multimathop
ad dist div dom meas sign supp .

\def\accorpa #1#2{\eqref{#1}--\eqref{#2}}
\def\Accorpa #1#2 #3 {\gdef #1{\eqref{#2}--\eqref{#3}}%
  \wlog{}\wlog{\string #1 -> #2 - #3}\wlog{}}


\def\tonde #1{\left(#1\right)}

\def\graffe #1{\mathopen\{#1\mathclose\}}

\def\<#1>{\mathopen\langle #1\mathclose\rangle}
\def\norma #1{\mathopen \| #1\mathclose \|}

\def\iot {\int_0^t}
\def\ioT {\int_0^T}
\def\intQt{\int_{Q_t}}

\def\intQ{\int_Q}
\def\iO{\int_\Omega}

\def\intQtau {\int_{Q_\tau}}

\def\dt{\partial_t}
\def\dn{\partial_{\bf n}}
\def\dtt{\partial_{tt}}

\def\checkmmode #1{\relax\ifmmode\hbox{#1}\else{#1}\fi}

\def\aeQ{\checkmmode{a.e.\ in~$Q$}}

\def\aat{\checkmmode{for a.a.~$t\in(0,T)$}}


\def\erre{{\mathbb{R}}}

\def\enne{{\mathbb{N}}}




\def\genspazio #1#2#3#4#5{#1^{#2}(#5,#4;#3)}
\def\spazio #1#2#3{\genspazio {#1}{#2}{#3}T0}

\def\L {\spazio L}
\def\H {\spazio H}
\def\W {\spazio W}

\def\C #1#2{C^{#1}([0,T];#2)}


\def\Lx #1{L^{#1}(\Omega)}
\def\Hx #1{H^{#1}(\Omega)}

\def\Ldue{\Lx 2}

\def\Huno{\Hx 1}
\def\Hdue{\Hx 2}



\let\theta\vartheta
\let\eps\varepsilon
\let\phi\varphi
\let\lam\lambda

\let\TeXchi\chi                         
\newbox\chibox
\setbox0 \hbox{\mathsurround0pt $\TeXchi$}
\setbox\chibox \hbox{\raise\dp0 \box 0 }
\def\chi{\copy\chibox}



\def\bQ{b_1}
\def\bO{b_2}
\def\bQh{b_3}
\def\bOh{\gamma}
\def\bz{b_0}
\def\bmas{b_4}
\def\btime{b_5}

\def\phQ{\phi_Q}

\def\phO{\phi_\Omega}


\def\sO{\s_\Omega}

\def\sQ{\s_Q}

\def\Uad{\calU_{\ad}}
\def\uopt{\overline u}
\def\tauopt{\overline{\tau}}

\def\Vp{V^*}

\def\normaV #1{\norma{#1}_V}

\let\hat\widehat

\def\Pi{\hat\pi}


\def\cd{c_\delta}
\def\s{\sigma}  
\def\m{\mu}	    
\def\ph{\phi}	
\def\a{\alpha}	
\def\b{\beta}	
\def\d{\delta}  
\def\et{\eta}   
\def\th{\theta} 
\def\r{\rho}    
\def\bph{\overline\ph}  
\def\bm{\overline\m}    
\def\bs{\overline\s}    
\def\ch{\chi}      
\def\J{{\cal J}}
\def\Jred{{\J}_{\rm red}}
\def\S{{\cal S}}
\def\I2 #1{\int_{Q_t}|{#1}|^2}
\def\IN2 #1{\int_{Q_t}|\nabla{#1}|^2}
\def\IO2 #1{\iO |{#1(t)}|^2}
\def\INO2 #1{\iO |\nabla{#1}(t)|^2}
\def\UR{{\cal U}_R}
\def\CP{{(CP)}}

\def\ts {\tau_*}
\def\x{{\bf x}}

\Begin{document}
\title{
		Penalisation of Long Treatment Time and Optimal Control of a Tumour Growth Model of Cahn--Hilliard
		Type with Singular Potential }
\author{}
\date{}
\maketitle

\Bcenter
\vskip-1cm
{\large\sc Andrea Signori$^{(1)}$}\\
{\normalsize e-mail: {\tt andrea.signori02@universitadipavia.it}}\\[.25cm]
$^{(1)}$
{\small Dipartimento di Matematica e Applicazioni, Universit\`a di Milano--Bicocca}\\
{\small via Cozzi 55, 20125 Milano, Italy}

\Ecenter
\Begin{abstract}\noindent
A distributed optimal control problem for a diffuse interface 
model, which physical context is that of tumour growth dynamics, is addressed.
The system we deal with comprises a Cahn--Hilliard equation
for the tumour fraction coupled with a reaction-diffusion for 
a nutrient species surrounding the tumourous cells.
The cost functional to be minimised possesses some objective
terms and it also penalises long treatments time, which may affect harm to the patients, 
and big aggregations of tumourous cells.
Hence, the optimisation problem leads to the optimal strategy which reduces the time exposure 
of the patient to the medication and at the same time allows the doctors to 
achieve suitable clinical goals.

\vskip3mm

\End{abstract}

\vskip3mm

\pagestyle{myheadings}
\newcommand\testopari{\sc Signori}
\newcommand\testodispari{\sc {}}
\markboth{\testodispari}{\testopari}

\salta
\finqui

\section{Introduction}
\label{SEC_INTRODUCTION}
\setcounter{equation}{0}
In the last decades, several developments have been obtained by scientists
in the field of tumour growth modelling. The key idea behind these models 
arises from realising that the tumour
tissue, as a special material, has to obey physical laws.
Hence, the modelling techniques originally developed for
engineering purposes can be adapted and exploited to 
derive mathematical models which
better emulate the evolution of tumours (see \cite{CL}).
The great advantages of mathematics are, among others, 
that of being able to foresee, make predictions, 
and capture information that does not interfere with the patient's health.
Moreover, mathematics
has the ability to select specific mechanisms we could be interested in.
Besides, let us also mentioned that
further understanding from the mathematical point of view
can also allow the doctors to tailor a personalised therapeutic pathway.

The diffuse interface model we are going to deal with reads as follows:
\begin{align}
   \a \dt \m + \dt \ph - \Delta \m &= P(\phi) (\sigma - \mu)
  \quad &&\hbox{in $\, Q:= \Omega \times (0,T)$},
  \label{EQ_prima}
  \\
   \mu &= \beta \dt \phi - \Delta \phi + F'(\phi)
  \label{EQ_seconda}
  \quad &&\hbox{in $\,Q$},
  \\
  \dt \sigma - \Delta \sigma  &= - P(\phi) (\sigma - \mu) + u
  \label{EQ_terza}
  \quad&& \hbox{in $\,Q$},
  \\
   \dn \m &= \dn \ph =\dn \s = 0
  \quad &&\hbox{on $\,\Sigma:= \partial \Omega \times (0,T)$},
  \label{EQ_BC}
  \\
   \m(0)&=\m_0,\, \ph(0)=\ph_0,\, \s(0)=\s_0
  \quad &&\hbox{in $\,\Omega,$}
  \label{EQ_IC}
\end{align}
\Accorpa\EQ EQ_prima EQ_IC
where $\a$ and $\b$ represent two positive relaxation parameters and $\Omega$ and $T>0$ denote
the spatial set in which the evolution takes place and the time horizon, respectively. 

This model constitutes a variation of the four-species thermodynamically
consistent model proposed
by Hawkins--Daruud et al. in \cite{HDZO} (see also \cite{HKNZ,CLLW,HDPZO,OHP}),
where the velocity contribution and chemotaxis 
effect are neglected and two relaxation terms $\a \dt \m$ and $\b\dt\ph$ in equation \eqref{EQ_prima} and \eqref{EQ_seconda} are included.
Let us noticing that equations \eqref{EQ_prima}--\eqref{EQ_seconda} comprises
of a viscous Cahn--Hilliard system for $(\ph,\m)$ (see, e.g., the review article \cite{Mir_CH} and the references therein for more details) with 
a non-standard source term $P(\phi) (\sigma - \mu)$ 
modelling the growth and death of cells.
Since the physical background of the above model has been extensively
described in \cite{FGR,CGH,CGRS_VAN,CGRS_ASY},
we just sketch the role covered by the occurring symbols.
In the above equations, the primary variables of the model are $\ph$, $\m$, and $\s$ denoting
in the order the difference in volume fractions between the tumour and healthy cells, the associated chemical potential, and the nutrient concentration of an unknown nutrient species (e.g., glucose, oxygen).
Typically $\ph$ ranges between 
$-1$ and $1$, where the two extremes represent the pure phases, i.e. the healthy case
and the tumourous case.
The function $P$ represents a source/sink term which
accounts for biological mechanisms such as proliferation.
Besides, as typical for phase field models, the function $F$ denotes a double-well potential whose classical examples are
the regular quartic potential and the singular logarithmic potential, 
which is more relevant for
applications. They are defined, in the order, by
\begin{align}
	& F_{reg}(r):= \tfrac 14 (r^2-1)^2 
	\quad \hbox{for } r \in \erre,
	\label{F_reg}
	\\ 
	& F_{log}(r):= 
	((1-r)\log(1-r)
	+(1+r)\log(1+r))
	- \lam r^2
	\quad \hbox{for } |r|<1,
	\label{F_log}
\end{align}
where in \eqref{F_log} $\lam$ stands for a positive constant large enough to avoid convexity.
The positive constants $\a$ and $\b$ 
can be seen as relaxation parameters.
The first one provides \eqref{EQ_prima} of a parabolic structure,
whereas the second term in the equation \eqref{EQ_seconda} is the classical viscous term of the Cahn--Hilliard equation. Let us also refer to \cite{S} and to \cite{CGH,CGRS_ASY,CGRS_VAN}
where these relaxations are incorporated in the model.
Lastly, the variable $u$ appearing in 
\eqref{EQAggterza} plays the role of control variable so that the system \EQ\ will be referred to as the {\it state} system in the sequel.

The well-posedness and long-time behavior of the above model (with $u=0$), 
in terms of the omega-limit set, have been addressed
in \cite{CGH} for a general class of double-well potentials, in the case $\a=\b>0$.
Next, in \cite{CGRS_VAN} and \cite{CGRS_ASY} Colli et al. discuss in which
sense the parameters $\a$ and $\b$ can be let to zero both separately and jointly
providing to specify the functional framework under which this can be done
depending on the asymptotic study under consideration.
Furthermore, we mention \cite{FGR}, where the above system (with $u=0$)
without any relaxation terms, i.e. the system \EQ\ with $\a=\b=0$, is studied. 
There, despite they restricted the analysis on
regular potentials with suitable polynomial growth, they keep the assumptions on the nonlinearity 
$P$ very general postulating for it a controlled polynomial growth.
Then, we refer to \cite{MRS}, where the authors investigate the 
existence of the global attractor for the dynamical system generated by \EQ~in the case $\a=\b=0$
(see, .e.g. \cite{MirZel} for further details on global attractors).
Furthermore, let us point out
\cite{FRL}, where a non-local model is taken into account
for the challenging case of singular potentials and degenerate mobilities.
As for the diffuse interface models including 
the velocity field effects by assuming
a Darcy's law or a Stokes--Brinkman's law, we refer to
\cite{WLFC,GLSS,DFRGM ,GARL_1,GARL_4,GAR, EGAR, FLRS,GARL_2, GARL_3}, 
where also further biological mechanisms such as chemotaxis and active transport are incorporated in the model. 
Lastly, we refer to \cite{Agosti,WZZ} and to the reference therein for some
numerical applications.

Before introducing the optimal control problem we are going to address, let us spend some words 
explaining how the cancer treatments are usually scheduled which
motivate some of the modelling considerations made below.
A typical medical treatment include surgery, chemotherapy, radiotherapy, and
immunotherapy and especially the last three therapies are particularly sensitive to the time 
exposition of the patient. 
Moreover, the therapy is divided into cycles consisting of a ``short" period of
treatment followed by a longer period of rest.
The goal of these therapies is usually to reduce the tumour mass
to achieve a reasonable stage which is compatible with surgery.
Besides, as time passes by the dispensed drug
starts to accumulate in the body bringing additional
waste items to be purified by kidneys and liver, and in the worst-case scenario, it may happen that 
after a long-time exposure the tumour cells became resistant to the medicament.
This is the reason which leads us to incorporate a long treatment time penalisation
in our optimisation problem.

Up to our knowledge,
the first contribution concerning an optimal control problem
governed by the system \EQ\ is \cite{CGRS_OPT}, where
\EQ\ was considered without any relaxation terms, i.e. with the choice $\a=\b=0$.
There, the authors proved the existence of a minimiser for the optimisation problem
and provide first-order necessary conditions for optimality 
in the framework of regular potentials exhibiting polynomial growth.
Then, we mention \cite{S}, where a similar optimisation problem is addressed for
the system \EQ\ with $\a,\b>0$ as the state system in
the case of singular, while regular, potentials allowing the
logarithmic potential to be included in the investigation. In this direction, the artificial relaxation terms
$\a\dt\m$ and $\b\dt\ph$ play a crucial role since, owing to their regularizing effect,
they allow proving a uniform separation principle for the phase variable 
which is a key property to
handle singular potentials. Next, the same author extends in \cite{S_DQ} the optimisation problem to the case of the double-obstacle potential 
by following the asymptotic scheme known in the literature as to {\it deep quench} limit.
Then, in \cite{S_a} and \cite{S_b} 
the author proves, by employing proper asymptotic strategies,
how the optimal control problem for the case
$\a,\b>0$ can be useful to solve the optimal control problems related to the state
system above in which $\a=0,\b>0$ and $\a>0,\b=0$
by letting the parameters $\a$ and $\b$ go to zero, respectively.
Besides, we are also aware of the recent work
\cite{CRW}, where, after discussing the long-time behavior of solutions, the authors
show that the optimal control problem \cite{CGRS_OPT} can be extended to the case in which
the cost functional also depends on time.
Referring to different models, we mention the contribution \cite{GARLR},
where an optimal treatment time has been performed for a slightly 
different state system of Cahn--Hilliard type, where the control appears in the first 
equation.
Moreover, we refer to \cite{SW}, where an optimal control problem
for the two-dimensional Cahn--Hilliard--Darcy system
with mass sources is addressed.
Lastly, we point out \cite{EK_ADV,EK}, where optimal control problems 
for the more involved
Cahn--Hilliard--Brinkman model, previously investigated by \cite{EGAR}, are addressed. 
For the interested reader, we also mention \cite{SM}, where a different type of 
control problem, known as {\it sliding mode} control, is performed for a
similar system.

In the spirit of \cite{GARLR} (and \cite{CRW}), we aim at generalizing the 
results established in \cite{S}
by introducing in the functional to be minimised a time penalisation. 
The optimisation problem considered in \cite{S} consists in 
minimising the following cost functional
\begin{align}
	 \non
	{\J}_{\rm old} (\ph, \s, u)  & = 
	\frac {{b_1}}2 \intQ |\ph - \phQ|^2
	+\frac {{b_2}}2  \iO|\ph(T)-\phO|^2
	+ \frac  {{b_3}}2 \intQ  |\s - \sQ|^2
	\\ & \quad\quad
	+ \frac  {{b_4}}2 \iO|\s(T)-\sO|^2
    + \frac  {{b_0}}2 \intQ |{u}|^2,
    \label{costfunct_old}
\end{align}
for some non-negative constants ${b_0},...,{b_4}$ and some
given target functions $\phQ,\sQ, \phO,\sO$ 
defined in proper functional spaces,
under the constrained that
the control $u$ belongs to the set of {\it admissible controls} $\Uad$ which is defined by
\Beq
    \label{Uad}
	\Uad := \graffe{u \in L^{\infty}(Q): u_* \leq u \leq u^*\ \aeQ},
\Eeq
where $ u_* $ and $ u^*$ denote some prescribed functions in $L^\infty(Q)$,
and such that the variables $\ph$ and $\s$ are solutions to the state system \EQ.
Thus, we aim at extending the above minimisation problem to a time-dependent cost functional 
by adding a free terminal time, which penalises long treatments time,
as well as an objective time to be approached.
Moreover, we also introduce in the cost functional an additional penalisation term
for large aggregation of tumour cells.
Namely, the time-dependent objective cost functional
we are going to minimise reads as 
\begin{align}
	 \non
	\J (\ph, \s, u, \tau)  & := 
	\frac \bQ 2 \intQtau |\ph - \phQ|^2
	+ \frac \bO 2 \iO |\ph(\tau)-\phO|^2
	+ \frac \bQh 2 \intQtau  |\s - \sQ|^2
	\\ & \quad\quad
	+\frac \bmas 2 \iO (1+\ph(\tau))
	+ \btime \tau
	+ \frac {b_6}2 |\tau - \ts|^2
    + \frac \bz 2 \intQ |{u}|^2,
    \label{costfunct}
\end{align}
where the symbols $b_0,...,b_6$ denote non-negative constants, while
$\phQ,\sQ,\phO$, and $\ts$ stand for the targets we want to approximate.
Here, let us point out the following comments
and differences with respect to the problem discussed in \cite{S}:
\begin{itemize}
\item[(i)] 
	Despite the last term in \eqref{costfunct}, the time integrals are performed between zero and
	an arbitrary $\tau \in [0,T]$,
	where $\tau$ models the treatment time of the cycle
	which the patient undergoes the clinical therapy,
	while $T$ may be regarded as the maximum amount of time 
	prescribed by some protocol.
	Let us claim that
	only minor changes are in order if one substitutes the term $\btime \tau$ in \eqref{costfunct} with a more general term like $\btime f(\tau)$, where 
	$f:[0,\infty)\to[0,\infty)$ is an increasing and continuously differentiable
	function.
\item[(ii)] The term $\frac {b_6}2 |\tau - \ts|^2$ forces the optimal time to
	be as close as possible to $\ts$ which stands for some target time to be reached.
	Notice that the sum of this latter (which is quadratic in $\tau$) and $\btime \tau$ (which is linear)
	gives still a convex contribution.
\item[(iii)]
	Minimizing the integral $\intQtau |\ph - \phQ|^2$ leads the phase variable $\ph$
	to be as close as possible, at the time $\tau$
	and in the sense of the $L^2$-norm, to the prescribed target $\phQ$. 
	In a similar fashion, it goes for the other variables. 
	Thus, $\phQ,\sQ,\phO$, should be chosen
	as stable configurations of the system or as some desirable configurations
	which are meaningful for surgery.
\item[(iv)] 
	The last term $ \intQ |{u}|^2$ penalises the 
	large values of the control variable designing the 
	side-effect that the dispensation of too many drugs to the patient
	may cause. 
\item[(v)]
	The term $\frac {1+\ph(\tau)} 2$ measures the size of the tumourous mass at
	the given time $\tau$. Hence, it penalises
	the strategies which do not shrink the tumour. Notice
	that the presence of $1$ in the numerator 
	is due to the fact that in the healthy case we have $\ph=-1$
	so that in that case the corresponding tumour mass is indeed zero.
\item[(vi)]
	The constants $b_0,...,b_6$
	can be chosen accordingly to the therapeutic goal we are interested in.
	\end{itemize}
Compared to \cite{GARLR} and \cite{CRW}, let us underline that we also include in the analysis a target 
time $\ts$ to be approximated as best as possible. This choice has the
advantage to produce a better characterisation of the optimality of the time
variable (cf. Theorem \ref{THM_Opt_cond_prima}). Notice that the choice $\ts=0$ is allowed.

To conclude the section, let us introduce some general facts concerning optimal control theory.
At first, let us note that since the well-posedness of system \EQ\ 
has already been established in \cite{S},
we are in a position to properly define the {\it control-to-state} operator
which assigns to a given control $u$ the corresponding solution to \EQ.
Namely, we have
\Beq
	\label{controltostate}
	\S:u \mapsto (\m,\ph,\s),
\Eeq
where $(\m,\ph,\s)$ stands for the unique solution to system \EQ~associated with the
control variable $u$.
This allows us to  suppress the variables $\ph$ and $\s$ in the cost functional 
$\J$ by expressing them as functions of $u$ leading to the corresponding {\it reduced}
cost functional which is defined as
\Beq
\label{J_red}
	\Jred(u,\tau):= \J (\S_2(u),\S_3(u),u,\tau),
\Eeq
where $\S_2(u),\S_3(u)$ denote the second and third components of $\S$, respectively.
Although the existence of a minimiser of the above problem can be deduced by following 
similar reasoning as in \cite{S}, the 
corresponding first-order necessary conditions for optimality 
present significant differences.
However, following classical arguments (see, e.g., \cite{Trol,Lions_OPT}), it is clear that
that the optimality of $(\overline{u},\overline{\tau})\in \Uad \times [0,T]$
can be characterised by employing the following variational inequalities
\Beq
	\label{optimal_formal}	
	\begin{cases}
	D_u \Jred(\overline{u},\overline{\tau})(v-\overline{u},\overline{\tau}) \geq 0 \quad \hbox{for every $v \in \Uad$},
	\\
	D_{\tau} \Jred(\overline{u},\overline{\tau})(\overline{u},s- \overline{\tau}) \geq 0 \quad \hbox{for every $s \in [0,T]$},
\end{cases}
\Eeq
where $D_{i}\Jred$, $i=\{u,\tau\}$, stand for the derivative 
of the reduced cost functional $\Jred$ with respect to the 
corresponding variable in a proper functional setting.
 
To get the necessary conditions for the time optimality,
the key argument is to show that 
the reduced cost functional is
\Frechet\ differentiable so that the abstract conditions \eqref{optimal_formal}
can be exploited.
In particular, the \Frechet\ differentiability with respect to the time variable requires higher order temporal regularity for the phase variable $\ph$.
In this direction, let us anticipate that the time derivative of the reduced cost functional 
will produce some terms involving
$\bph(\overline{\tau})$ and $ \dt \bph(\overline{\tau})$ 
(cf. Theorem~\ref{THM_frechet_time}), where
$(\uopt, \overline{\tau})$ stands for some optimal pair and $(\bm,\bph,\bs)$ for 
the corresponding state.
A sufficient condition which gives meaning to the above terms is
$\bph \in \H2 {\Lx2}$ due to the continuous injections 
of $\H2 {\Lx2} $ in $\C1 {\Lx2}$ so that the pointwise
terms $\bph(\overline{\tau})$ and 
$ \dt \bph(\overline{\tau})$ are meaningful, at least in ${\Lx2}$. 
In this regards, let us also mention that we can not consider in the cost functional \eqref{costfunct}
any contribution involving $\iO |\s(\tau)-\sO|^2$,
where $\sO$ models some target function.
The reason is that if such term is present, in the time derivative of $\J_{\rm red}$
it will appear the pointwise term $\dt\bs(\overline{\tau})$
which, in turn, will require to show $\bs \in \H2 {\Lx2}$. However,
the nutrient equation \eqref{EQAggterza} contains the control variable $u$ 
so that, to get the mentioned regularity for $\s$, 
we would be forced to assume the control variable to be sufficiently regular in time, say 
$u\in\H1 {\Lx2}$, which is not significant for the applications.
For the same reason, we are considering the term $\intQ |u|^2$ in \eqref{costfunct}, instead of
$\int_{Q_\tau} |u|^2$, to avoid assuming any temporal regularity for $u$.
However, in Section \ref{SEC_GENERALIZING_COST_FZ}, by using the relaxation arguments employed by Garcke et al. in \cite{GARLR},
we also show that it is possible to include in the analysis an objective control for the nutrient variable at time $\tau$.

Hence, introducing the space of admissible states and controls by
\Beq
	{\cal A}_{ad}:= \biggl\{\bigl(\ph,\s,u, \tau\bigr) : (u,\tau) \in \Uad \times [0,T], 
	\hbox{ such that } (\ph,\s)=\bigl(\S_2(u),\S_3(u)\bigr) \biggr\},
	\label{AD}
\Eeq
we can summarise the  minimisation problem we are going to address as:
\begin{align*}
	\CP & 
	\quad 
	\inf_{(\ph,\s,u, \tau) \in {\cal A}_{ad}}
	\J(\ph,\s,u, \tau).
\end{align*}

The rest of the paper is outlined as follows: in the next section, we set our conventions, present
the assumptions and state our results.
The existence of a minimiser for the optimisation problem $\CP$ and 
the corresponding first-order necessary conditions for optimality have been addressed in 
Section \ref{SEC_EXISTENCE:OPTIMAL_CONTROL}. Next, in Section \ref{SEC_GENERALIZING_COST_FZ}, we point out some
possible generalisations of the work via a relaxation argument.

\section{Mathematical Setting}
\label{SEC_MATH_SETT}
\setcounter{equation}{0}
Throughout the paper, we assume $\Omega\subset\erre^3$ to be a smooth bounded 
domain with boundary $\Gamma$, and $T>0$ is a fixed final time.
For every $t\in[0,T]$, we employ the classical notation 
\begin{align*}
	Q_{t}& :=\Omega \times (0,t), \quad \Sigma_{t}:=\Gamma \times (0,t)
	\quad \hbox{for every $t \in (0,T]$,}
\end{align*}
and
\begin{align*}
	Q:=Q_{T}, \quad \Sigma:=\Sigma_{T}.
\end{align*}
For an arbitrary Banach space $X$, we use
$\norma{\cdot}_{X}$ to denote its norm, $X^*$ 
for its topological dual, and ${}_{X^*}\<\cdot, \cdot >_X$ 
for the duality product between $X^*$ and $X$.
Meanwhile, for every $p \in [1,+\infty]$,
we simply write $\norma{\cdot}_{p}$ to indicate the usual norm of the Sobolev spaces $L^p(\Omega)$.
Besides, it turns out to be convenient to set 
\Beq
	\non
	H:= \Ldue, \quad V:= \Huno, \quad W:=\graffe{v \in \Hdue : \dn v = 0 \hbox{ on } \Gamma},
\Eeq
equipped with their standard norms,
where $\dn$ stands for the outward normal derivative 
of $\Gamma$.
Under these assumptions, it follows that the injections
$V \hookrightarrow H \cong H^* \hookrightarrow V^*$ 
are both continuous and dense which entails that
$(V,H,\Vp)$ forms a Hilbert triple so that 
we have the following identification
\Beq
	\non
	{}_{\Vp}\< u,v >_V = \iO u v \quad \hbox{for every $u \in H$, $v \in V$.}
\Eeq
As far as the general assumptions are concerned, we postulate that
\begin{align}
	& 
	\a, \b \, \hbox{ are positive constants}.
	\label{ab}
	\\ &
	\bz, \bQ, b_2, b_3, b_4 ,b_5, b_6 \, \hbox{ are non-negative constants, but not all zero}.
	\label{constants}
	\\ &
	\phQ, \sQ: Q \to \erre, \phO :\Omega \to \erre, \hbox{ and }
	\phQ, \sQ \in L^2(Q), \phO \in \Lx2.
	\label{targets}
	\\
	& u_*,u^* \in L^\infty(Q) \hbox{ with } u_*\leq u^* \,\aeQ.
	\\
	&  \ts \in [0,T].
	\\ 	\non &
	P \in C^2(\erre) \hbox{ is non-negative, bounded with bounded derivative and}
	\\
	& \hbox{Lipschitz continuous}.
	\label{P}
	\\ &
	\ph_0 \in W, \m_0 \in V\cap\Lx\infty, \s_0 \in V.
	\label{initial_data}	
	\\ &
	\label{F_L1}
	 F(\ph_0)\in \Lx1.
\end{align}
Moreover, let us postulate that the control-box $\Uad$ is defined by \eqref{Uad},
so that $\Uad$ is a closed and convex subset of $L^2(Q)$.
On the other hand, it will be sometimes convenient to work with  
an open set. Hence, let us define the open superset $\UR$ as follows
\begin{align}
	& \non
	\UR \subset L^2(Q) \hbox{ is a non-empty, bounded and open set containing } \Uad
	\\ \non &
	 \hbox{such that } \norma u_2 \leq R \hbox{ for all }
	u \in \UR.
\end{align}
As for the nonlinear double-well potential $F$, we require that
\Beq
	\label{defF}
	F:\erre \to [0, +\infty), \quad  \hbox{with} \quad  
	F:= \hat{B} + \hat{\pi},
\Eeq
where
\begin{align}
	& 
	\hat{B}:\erre \to [0, +\infty] \ \ \hbox{is convex, and lower semicontinuous, }	
	\label{Bhat}
	\hbox{with {$\hat{B}(0)=0$}}.
	\\ &
	\hat{\pi} \in C^3(\erre) \,\, \hbox{and} \,\, \pi:=\hat{\pi}' 
	\hbox{ is \Lip~continuous.}
	\label{pihat}
\end{align}
Under these assumptions it is well-known that $ B:=\partial\hat{B}$ is a maximal and monotone 
graph $B \subseteq \erre \times \erre$
(see, e.g., \cite[Ex. 2.3.4, p. 25]{BRZ}) whose domain we indicate by $D(B)$.
Furthermore, we assume that $F$ is a smooth function when restricted to its domain
by assuming that
\begin{align}
	 \non
	\ D(B)&=(r_-,r_+), \hbox{ with } -\infty \leq r_- <0 < r_+ \leq + \infty.
	\\
	\label{F_dominio}
	 F_{| _{D(B)}} &\in C^3(r_-,r_+), \hbox{ and } \lim\limits_{r \to r_{\pm}} F'(r) = \pm \infty.
\end{align}
It is worth noting that both the regular potential \eqref{F_reg}
and the logarithmic potential \eqref{F_log}
do fit the above assumptions.
Moreover, we additionally require that the initial datum $\ph_0$ verifies
\begin{align}
	& \label{phzeropiccolo}
	r_- < \inf \ph_0 \leq \sup \ph_0 < r_+,
\end{align}
which from the physical viewpoint means that $\ph_0$ does not contain any region with pure phases.
Notice that the above condition combined with \eqref{initial_data} entails that
\begin{align}
	\label{dtphzeroinH}
	\tfrac 1\b \tonde{\m_0 + \Delta \ph_0 - B(\ph_0) - \pi(\ph_0)} \in H.
\end{align}
\Accorpa\tutteleipotesi ab dtphzeroinH
The mathematical assumptions required so far are more or less the same assumed in \cite{S}. 
However, as already mentioned, the first-order necessary conditions
for optimality
that we will point out later will demand higher order temporal regularity for the phase variable.
To give meaning to all the appearing pointwise terms we
replace \eqref{targets} and \eqref{dtphzeroinH} by
\begin{align}
	& \label{targets_reg}
	\phQ, \sQ \in \H1 H,
	\\ &
	\label{dtphzeroinV}
	\tfrac 1\b \tonde{\m_0 + \Delta \ph_0 - B(\ph_0) - \pi(\ph_0)} \in {V},
\end{align}
\Accorpa\Regdata targets_reg dtphzeroinV
respectively. Let us also point out that condition \eqref{dtphzeroinV} easily follows 
once the initial data, in addition to \eqref{initial_data}, fulfils $\ph_0\in \Hx3$.

Before moving on, let us recall some \wk~results which we will apply later on.
First, we often owe to the standard Sobolev continuous embedding
\begin{align}
	\label{VinL6}
	\Huno \hookrightarrow L^q(\Omega),
	\quad \hbox{for every $q\in [1,6]$,}
\end{align}
which is also compact for every $q \in [1,6)$.
Moreover, we recall the Young inequality
\begin{align}
  ab \leq \delta a^2 + \frac 1 {4\delta} \, b^2,
  \quad \hbox{for every $a,b\geq 0$ and $\delta>0$}.
  \label{young}
\end{align}
Lastly, we convey to use the symbol small-case $c$ for every constant
which only depend on structural data of the problem such as
the final time $T$, $\Omega$, $R$, the shape of the 
nonlinearities, the norms of the involved functions, and possibly $\a$ and $\b$. 
On the other hand, we devote the capital letters to designate 
some specific constants.
\section{The Control Problem}
\label{SEC_EXISTENCE:OPTIMAL_CONTROL}
\setcounter{equation}{0}
\subsection{The State System}
To begin with, let us recall the well-posedness result for the sytem \EQ\ obtained in \cite{S}.
\Bthm[{{\cite[Thms.~2.1,~2.2, and~2.3]{S}}}]
\label{THM_Wp_state_system}
Suppose that \tutteleipotesi\ hold and let $u \in \UR$. 
Then, the state system \EQ\ admits a unique solution $(\m,\ph,\s)$ satisfying
\begin{align}
	\ph &\in W^{1,\infty}(0,T;H) \cap \H1 V \cap \L\infty W \subset \C0 {C^0({\overline{\Omega}})},
	\label{reg1}
	\\ 
	\m, \s &\in \H1 H  \cap \L\infty V \cap \L2 W \subset C^0([0,T];V),
	\label{reg2}
	\\ 
	\m &\in L^{\infty}(Q).
	\label{reg3}
\end{align}
\Accorpa \reg reg1 reg3
Moreover, there exists a positive constant $C_1$, which depends on $R$, $\a$, $\b$, 
and on the data of the system, such that
\begin{align}
 	& \non 
 	\norma{\ph}_{W^{1,\infty}(0,T;H) \cap \H1 V \cap \L\infty W}
	+ \norma{\m}_{\H1 H  \cap \L\infty V \cap \L2 W \cap L^\infty(Q)}
	\\ & \quad
	+ \norma{\s}_{\H1 H \cap \L\infty V \cap \L2 W}
	\leq C_1.
	\label{regstima}
\end{align}
In addition, it holds the \socal~uniform separation property. 
Namely, there exists a 
compact subset $K\subset (r_-,r_+)=D(B)$ such that
\Beq
	\non
	\ph(\x, t) \in K \quad \hbox{for all } (\x,t) \in Q.
\Eeq
Furthermore, the following estimate
\Beq
	\label{stimasep}
	\norma{\ph}_{C^0(\overline{Q})}
	+ \max_{0 \leq i \leq 3} \,\norma{F^{(i)}(\ph)}_{L^\infty(Q)}
	+ \max_{0 \leq j \leq 2} \,\norma{P^{(j)}(\ph)}_{L^\infty(Q)}
	\leq C_2
\Eeq
is satisfied for a positive constant $C_2$ which depends only on $R$, $\a$, $\b$, 
$K$ and on the data of the system.
\Ethm
The well-posedness of the state system \EQ\ established by the the above theorem allow us to define
the control-to-state operator $\S$ as the map which assigns to every control
$u$ the corresponding solution $(\m,\ph,\s)$ to system \EQ.

We are now ready to present the first novelty of the work regarding 
improved regularity results for the solutions to \EQ\ obtained in the above theorem
that will be used later on to investigate the optimal control problem $\CP$.
\Bthm
\label{PROP_timereg}
Suppose that \accorpa{ab}{phzeropiccolo} and \eqref{dtphzeroinV} hold and let $u\in \UR$. 
Then, the unique solution $(\m,\ph,\s)$ to
\EQ\ obtained from Theorem \ref{THM_Wp_state_system}, in addition to \accorpa{reg1}{reg3}, enjoys the following regularity
\begin{align}
	\label{prop_time}
		\ph &\in \W{1,\infty} V \cap \H2 H \subset \C1 H\cap \C0 {\Hx2}.
\end{align}
Moreover, there exists a positive constant $C_{3}$ depending
on $R$, $\a$, $\b$, and on the data of the system, such that
\Beq
	\label{stima_new}
	\norma{\ph}_{\W{1,\infty} V \cap \H2 H} 
	\leq
	C_3.
\Eeq
\Ethm

\Bdim
For the sake of simplicity, we perform only formal a priori estimates
which can be carried out rigorously within an approximation scheme
such as the Faedo--Galerkin scheme. 
 
We differentiate \eqref{EQ_seconda} with respect to time, 
multiply the obtained equation by $\dtt \ph$, and integrate over $Q_t$ 
and by parts to obtain that
\begin{align}
	& \non
	\b \I2 {\dtt \ph}
	+ \frac 12 \IO2 {\nabla \dt \ph}
	=
	\frac 12 \iO |\nabla\dt\ph(0)|^2
	- \intQt F''(\ph)\dt \ph \, \dtt\ph
	+ \intQt \dt \m \, \dtt\ph,
\end{align}
where we denote the integrals on the \rhs\ by $I_1,I_2$ and $I_3,$ respectively.
The terms on the \lhs~are non-negative, whereas $I_2$ and $I_3$ can 
be dealt by means of the Young inequality, along with the estimate \eqref{stimasep} which is verified by the solution $\ph$. 
Namely, we have that
\begin{align}
	|I_2|+|I_3|
	& \non
	\leq
	\frac \b2 \I2 {\dtt \ph}
	+ c \intQt (|\dt \m|^2 + |\dt \ph|^2).
\end{align}
Moreover, by taking $t=0$ in \eqref{EQ_seconda} and using the assumption \eqref{dtphzeroinV}
we readily infer that
\begin{align}
	\non
	|I_1|\leq c.
\end{align}
Therefore, owing to the estimate \eqref{regstima}, we deduce that there exists a positive constant $c$ such that
\Beq
	\non
	\norma{\dtt \ph}_{\L2 H}
	+\norma{\nabla \dt \ph}_{\L\infty H}
	\leq
	c.
\Eeq
Next, let us recall the \wk~embedding of $\H1 H $ in $\C0 H$
which entails also that $\dt\ph \in \C0 H$.
Lastly, comparison in equation \eqref{EQ_seconda} produces
\begin{align*}
 	\Delta \ph \in \C0 H,
\end{align*}
so that, upon invoking the elliptic regularity theory, the proof is concluded.
\Edim

\Bthm
\label{PROP_timeregdue}
Suppose that the assumptions of Theorem \ref{PROP_timereg} 
are verified. Moreover, let $\m_0\in \Hx3$.
Then, the unique solution $(\m,\ph,\s)$ to
\EQ\ obtained from Theorem \ref{THM_Wp_state_system}, in addition to \accorpa{reg1}{reg3}
and \eqref{prop_time} satisfies
\begin{align}
		\m &\in \W{1,\infty} V \cap \H2 H \subset \C1 H \cap \C0 {\Hx2}.
\end{align}
Moreover, there exists a positive constant $C_{4}$ depending
on $R$, $\a$, $\b$, and on the data of the system, such that
\Beq
	\label{stima_new_due}
	\norma{\m}_{\W{1,\infty} V \cap \H2 H} 
	\leq
	C_4.
\Eeq
\Ethm
\Bdim
As before, we proceed formally and let us claim that the proof can be carried out in a rigorous fashion by employing a Galerkin scheme.

We differentiate \eqref{EQ_prima} with respect to time, test
the obtained equation by $\dtt \m$, and integrate over time and by parts to obtain that
\begin{align*}
	& \a \I2 {\dtt \m}
	+ \frac 12 \IO2 {\nabla \dt \m}
	=
	\frac 12 \iO |\nabla \dt \m(0)|^2
	+\intQt P'(\ph)(\s-\m)\dt\ph \dtt \m
	\\ & \quad
	+ \intQt P(\ph)(\dt\s-\dt\m)\dtt \m
	- \intQt \dtt\ph \dtt \m,
\end{align*}
where we indicate by $I_1,...,I_4$ the integrals on the \rhs, in this order.
The first term can be bounded by combining the assumption
\eqref{initial_data} with \eqref{dtphzeroinV}
and the additional requirement on the initial datum $\m_0 \in \Hx3$.
In fact, evaluating equation \eqref{EQ_prima} at $t=0$, and then
equation \eqref{EQ_seconda} at $t=0$ lead us to realise that
\begin{align*}
	\dt\m(0)=	
	\frac 1\a 
	\biggl(
	\frac 1\b \bigl( -\m_0 - \Delta \ph_0 + F'(\ph_0)\bigr)	
	+ \Delta \m_0 + P(\ph_0)(\s_0-\m_0)
	\biggr) \in V
\end{align*}
so that
\begin{align*}
	|I_1| \leq c.
\end{align*}
Then, using the previous estimates \accorpa{regstima}{stimasep} and \eqref{stima_new}, \Holder's inequality, the continuous embedding
$V\subset \Lx4$, and the boundedness of $P'$, we have
\begin{align*}
	|I_2|
	& \leq
	c \iot  (\norma{\s}_{4} + \norma{\m}_{4} )\norma{\dt\ph}_{4}\norma{\dtt\m}_{2}
	\\ & \leq 
	\d \I2 {\dtt\m}
	+ \cd \intQt (\norma{\s}_{V}^2 + \norma{\m}_{V}^2 )\norma{\dt\ph}_{V}^2
	\\ & \leq 
	\d \I2 {\dtt\m}+\cd,
\end{align*}
for a positive $\d$ yet to be determined,
where in the last line we also invoke the fact that, due to \accorpa{reg1}{reg3}
and to \eqref{stima_new}, we have that $\norma{\s}_{V},\norma{\m}_{V}$
and $\norma{\dt\ph}_{V}$ belong to $ L^\infty(0,T)$.
Using the Young inequality, the boundedness of $P$, and \eqref{stimasep},
we infer that
\begin{align*}
	|I_3|+|I_4|
	& \leq
	\d \I2 {\dtt\m}
	+ \cd \intQt (|\dtt \ph|+|\dt \s|^2+|\dt\m|^2),
\end{align*}
for a positive $\d$ yet to be determined.
Hence, adjusting $\d\in(0,1)$ small enough 
and accounting for the above estimates, we deduce that
\begin{align*}
	\norma{\dtt \m}_{\L2 H}
	+\norma{\nabla \dt \m}_{\L\infty H}
	\leq
	c.
\end{align*}
Arguing as above, we easily infer that $\dt\m \in \C0 H$
and then, by comparison in equation \eqref{EQ_prima} also that
\begin{align*}
	\Delta \m \in \C0 H
\end{align*}
which complete the proof upon invoking the elliptic regularity theory.
\Edim

To conclude the section, let us recall the continuous dependence result for \EQ\ 
obtained in \cite[Thms.~2.2~and~2.3]{S}.
\Bthm
\label{THM_Lip_conttroltostate}
Assume that \tutteleipotesi~are fulfilled. 
Moreover, for $i=1,2$, let $u_i \in \UR$ and $(\m_i,\ph_i,\s_i)$ be
the corresponding solution to \EQ\ obtained from Theorem \ref{THM_Wp_state_system}.
Then, there exists a positive constant $C_5$,
which depends only on $R$, $\a$ and $\b$, and on the data of the system such that
\begin{align}
	\non &
	\norma{\a(\m_1-\m_2) + (\ph_1-\ph_2) + (\s_1-\s_2)}_{\L\infty\Vp}
	+ \norma{\m_1-\m_2}_{\L\infty H \cap \L2 V}
	\\ \non & \qquad
	+\norma{\ph_1-\ph_2}_{\H1 H \cap \L\infty V}
	+ \, \norma {\s_1-\s_2} _{\L\infty H \cap \L2 V}
	\\ & \quad
	\leq C_5 \norma {u_1-u_2} _{\L2 H}.
\end{align}
\Ethm
 This latter can be equivalently 
interpreted as a \Lip~continuity property for the control-to-state operator $\S$
between suitable Banach spaces.

\subsection{Existence of a Minimiser}
\label{SEC_EX_MIN}
Here, we prove the existence of a minimiser for the optimisation problem $\CP$.
\Bthm
\label{THM_Existence_optimal_controls}
Assume that \tutteleipotesi~are fulfilled. Then, the optimal control problem
$(CP)$ admits at least a minimiser. Namely, there exists some
$(\bph,\bs,\overline{u},\overline{\tau})\in{\cal A}_{ad}$ such that
\Beq
	\non
	\J (\bph,\bs,\overline{u},\overline{\tau}) 
	= \inf_{(\ph,\s,u,\tau) \in {\cal A}_{ad}} \J(\ph,\s,u,\tau).
\Eeq
The control variable $\uopt$ will be referred to as optimal control,
whereas $\overline \tau$ and $(\bm,\bph,\bs)$ will be referred to as
optimal time and optimal state, respectively.
\Ethm
\Bdim
The proof can be easily carried out by means of the \wk~direct
method of calculus of variations, e.g. by retracing the proof of \cite[Thm. 2.6]{S}.
To begin with, let us check that the cost functional $\J$ is bounded from below.
Using the bounds for the phase variable $\ph$ expressed by \accorpa{regstima}{stimasep}, we
infer that
\begin{align}
	& \non
	\J (\ph, \s, u, \tau)   
	\geq  
	\frac \bmas 2 \iO \ph(\tau)
	\geq 
	- \frac \bmas 2 \norma{\ph}_{C^0(\overline Q)}
	\geq  - \frac { \bmas} 2 C_2
	> -\infty.
\end{align}
Next, we pick a minimizing sequence
$\graffe{(u_n,\tau_n)}_n$ of elements of $\Uad \times(0,T)$
with the sequence of the corresponding solutions $\{(\m_n,\ph_n,\s_n)\}_n$ 
to \EQ\ obtained from Theorem \ref{THM_Wp_state_system}.
Namely, we have that
\Beq
	\non
	\lim_{n\to\infty}\J(\ph_n,\s_n, u_n,\tau_n) 
	= \inf_{(\ph,\s,u, \tau) \in {\cal A}_{ad}}	\J(\ph,\s,u, \tau) > - \infty.
\Eeq
On the other hand, for every $n\in\enne$, the bounds provided 
by estimate \eqref{regstima} hold independentely of $n$.
Therefore, using classical weak and weak-star compactness results we infer that,
up to a not relabeled subsequence, there exist 
some $\overline{u} \in \Uad$ and a triple $(\bm,\bph,\bs)$ such that,
as $n\to \infty$,
\begin{align}
	& \non
	u_n \to \overline{u} \ \ \hbox{weakly star in } L^\infty(Q),
	\\ \non &
	\m_n \to \bm  \ \ \hbox{weakly star in } \H1 H  \cap \L\infty V \cap \L2 W \cap L^\infty(Q),
	\\ \non & 
	\ph_n \to \bph \ \ \hbox{weakly star in } W^{1,\infty}(0,T;H) \cap \H1 V \cap \L\infty W,
	\\ \non &
	\s_n \to \bs \ \ \hbox{weakly star in } \H1 H \cap \L\infty V \cap \L2 W.
\end{align}
Moreover, standard compactness arguments (see, e.g., \cite[Sec.~8, Cor.~4]{Simon}) 
yield that, possibily up to a not relabeled subsequence, as $n\to\infty$,
\Beq
	\label{stron_ph}
	\ph_n \to \bph \ \ \hbox{strongly in } C^0 (\overline{Q})
\Eeq
and also that there exists $\overline{\tau}\in [0,T]$ such that, as $n\to\infty$,
\Beq
	\label{tau_n}
	\tau_n \to \overline{\tau}.
\Eeq
Next, by using estimate \eqref{stimasep}, the strong convergence \eqref{stron_ph},
and the properties of $F$ and $P$, we realise that
\Beq
	\non
	F'(\ph_n) \to F'(\bph), \ \ P(\ph_n) \to P(\bph) \quad \hbox{ strongly in $C^0 (\overline{Q})$}.
\Eeq
Thus, it is then a standard matter to pass to the limit, as $n\to\infty$, in 
the variational formulation of \EQ\ written for $(\m_n,\ph_n,\s_n)$ and conclude that
$(\bm,\bph,\bs) = \S (\overline{u}) $.
Moreover, \eqref{tau_n} ensures that, as $n\to\infty$,
\Beq
	\label{chi_tau_n}
	\chi_{[0,\tau_n]} (t)\to \chi_{[0,\overline{\tau}]}(t) 
	\quad \hbox{$\aat$}.
\Eeq
Hence, we claim that the limit $(\bph,\bs,\overline{u},\overline{\tau})$ 
is indeed the minimiser we are looking for. 
Before showing how to pass to the limit term by term,
let us point out that
\Beq
	\non
	\int_{Q_{\tau_n}} |\cdot|^2
	=\int_0^{\tau_n} \norma{\cdot}_2^2
	= \int_0^T \norma{\cdot}_2^2 \, \chi_{[0,\tau_n]} .
\Eeq
As a matter of fact, it easily follows from the above convergences that
$\ph_n-\phQ \to \bph-\phQ$ strongly in $\L2 H$ so that, as $n\to \infty$,
\Beq
	\label{conv_existence_uno}
	\int_{Q_{\tau_n}} |\ph_n-\phQ|^2
	\to
	\int_{Q_{\overline{\tau}}} |\bph-\phQ|^2.
\Eeq
In fact, we have that
\begin{align*}
	& \ioT \biggl( \norma{\ph_n-\phQ}^2_2 \, \ch_{[0,\tau_n]} 
			-  \norma{\bph-\phQ}^2_2 \, \ch_{[0,\tauopt]} \biggr)
	\\ & \quad \leq 
	\ioT \norma{\ph_n-\phQ}^2_2 \bigl( \ch_{[0,\tau_n]}-\ch_{[0,\tauopt]} \bigr)
	+ \ch_{[0,\tauopt]} \ioT \biggl( \norma{\ph_n-\phQ}^2_2
			-  \norma{\bph-\phQ}^2_2 \biggr),
\end{align*}
where both the terms on the \rhs~go to zero
by combining the Lebesgue convergence theorem with
the pointwise convergence \eqref{chi_tau_n} and the strong convergence of $\ph_n-\phQ $.
Next, let us claim that the second term of the cost functional verifies that, as $n\to \infty$,
\Beq
	\label{primanew}
	\iO |\ph_n({\tau}_n)-\phO|^2
	\to
	\iO |\bph({\overline{\tau}})-\phO|^2.
\Eeq
In fact, it holds that
\begin{align*}
	& \non
	\Big|\iO |\ph_n({\tau}_n)-\phO|^2 -\iO |\bph({\overline{\tau}})-\phO|^2\Big|
	\leq \norma{\ph_n({\tau}_n)+\bph({\overline{\tau}})-2\phO}_2
	\norma{\ph_n({\tau}_n)-\bph({\overline{\tau}})}_2.
\end{align*}
Moreover, the convergences \accorpa{stron_ph}{tau_n}, along with
the triangular inequality and the fundamental theorem of calculus, 
allow us to handle the last term as 
\begin{align}	
	\non
	\norma{\ph_n({\tau}_n)-\bph({\overline{\tau}})}_2
	& \leq 
	\norma{\ph_n({\tau}_n)-\ph_n({\overline{\tau}})}_2
	+ \norma{\ph_n({\overline{\tau}})-\bph({\overline{\tau}})}_2
	\\ \non
	& \leq 
	|\tau_n - \overline{\tau}|^{\frac 12}
	\biggl(\int_{\overline{\tau}}^{{\tau}_n} \norma{\dt\ph_n}^2_{2}\biggr)^{\frac 12}
	+ \norma{\ph_n({\overline{\tau}})-\bph({\overline{\tau}})}_2
	\\ \non
	&\leq  
	|\tau_n - \overline{\tau}|^{\frac 12} \norma{\dt\ph_n}_{\L2 H}
	+ \norma{\ph_n({\overline{\tau}})-\bph({\overline{\tau}})}_2.
\end{align}
Notice that the first term on the \rhs~vanishes accounting for the bound 
\eqref{regstima} and for the convergence
\eqref{tau_n}. Meanwhile, the second term goes to zero due to the strong
convergence \eqref{stron_ph} so that \eqref{primanew} follows.
Furthermore, owing to \eqref{tau_n}, we easily infer that, as $n\to \infty$,
\begin{align*}
	|\tau_n-\ts|^2 \to |\tauopt -\ts|^2 .
\end{align*}
The remaining terms can be handled arguing in a similar fashion. 
Lastly, the weak sequential lower semicontinuity of $\J$, entails that 
\Beq
	\non
	\J (\bph,\bs,\overline{u},\overline{\tau}) 
	\leq 
	\liminf_{n\to\infty}\J(\ph_n,\s_n, u_n,\tau_n) 
	= \inf_{(\ph,\s,u, \tau) \in {\cal A}_{ad}}	\J(\ph,\s,u, \tau),
\Eeq
so that $(\bph,\bs,\overline{u},\overline{\tau})$ is indeed a minimiser
for $\CP$, as we claimed.
\Edim

\subsection{The Linearised System}
Once the existence of minimisers has been obtained, we aim at pointing out some first-order
necessary conditions for optimality by exploiting the theoretical conditions 
\eqref{optimal_formal}. In this direction, we first show the operator $\S$ is 
\Frechet~differentiable between suitable Banach spaces and then use the chain rule and the definition of
the reduced cost functional to develop the abstract variational inequalities \eqref{optimal_formal}.

The first step consists in investigating the linearised system of \EQ. 
For a fixed control $\overline{u}\in\UR$ with the corresponding state
$(\bm,\bph,\bs)$, and for an arbitrary $h\in L^2(Q)$ the linearised system to \EQ~reads as 
\begin{align}
   \a \dt \et + \dt \th - \Delta \et &= P'(\bph) (\bs - \bm)\th + P(\bph)(\r - \et)
  \quad &&\hbox{in $ Q$},
  \label{EQLinprima}
  \\
   \et &= \b \dt \th - \Delta \th + F''(\bph)\th
  \label{EQLinseconda}
  \quad &&\hbox{in $Q$},
  \\
   \dt \r - \Delta \r &= -P'(\bph) (\bs - \bm)\th - P(\bph)(\r - \et) + h
  \label{EQLinterza}
  \quad &&\hbox{in $Q$},
  \\
  \dn \r  &= \dn \th = \dn \et = 0
  \quad &&\hbox{on $\Sigma$},
  \label{BCEQLin}
  \\
   \r(0) &= \th(0) = \et (0) = 0
  \quad &&\hbox{in $\Omega$}.
  \label{ICEQLin}
\end{align}
\Accorpa\EQLin EQLinprima ICEQLin
The expectation is as follows: for every $h\in L^2(Q)$,
provided to find the proper Banach space ${\cal Y}$,
the operator $\S$ is \Frechet\ differentiable in ${\cal Y}$ and its
directional derivative along $h$ is given by the corresponding solution 
to the linearised system above, i.e. 
$D\S(\overline{u}) h = (\et, \th, \r)$. 
Since the linearised system is independent of the choice of the cost functional,
but only depends on the state system \EQ,
the well-posedness for the above system directly follows from \cite{S}.
\Bthm[{\cite[Thm.~2.4]{S}}]
\label{THM_Wp_linear}
Assume that \tutteleipotesi~are fulfilled. Then, for every h $\in L^2(Q)$, the linearised
system \EQLin~admits a unique solution $(\et,\th,\r)$ satisfying
\begin{align}
	\et, \th, \r \in \H1 H \cap \L\infty V \cap \L2 W \subset \C0 V.
	\label{reglin}
\end{align}
In addition, there exists a positive constant $C_6$, which depends on the data 
of the system, and possibly on $\a$ and $\b$, such that
\begin{align}
 	\non &
 	\norma{\et}_{\H1 H \cap \L\infty V \cap \L2 W}
	+ \norma{\th}_{\H1 H \cap \L\infty V \cap \L2 W}
	\\ & \non \quad
	+ \norma{\r}_{\H1 H \cap \L\infty V \cap \L2 W}
	\leq C_6.
\end{align}
\Ethm
Now, we can rigorously formulate our expectation concerning the 
\Frechet~differentiability of the map $\S$.
\Bthm[{\cite[Thm.~2.5]{S}}]
\label{THM_Frechet_space}
Assume that \tutteleipotesi~are satisfied and let $\uopt$
and $(\bm,\bph,\bs)$ be an optimal control for $\CP$ with the corresponding state. 
Then, the control-to-state operator $\S$ is 
\Frechet~differentiable at $\uopt$ as a mapping from $\UR$ into the space
$\cal Y$, where 
\Beq
	\label{statespaceFre}
	{\cal Y}:= \bigl( \H1 H \cap \L\infty V \cap \L2 W \bigr)^3.
\Eeq
Moreover, for any $\overline{u} \in \UR$, the \Frechet~derivative 
$D\S(\overline{u})$ is a linear and continuous operator from $L^2(Q)$ to $\calY$
such that 
\begin{align*}
	D\S(\overline{u}) h = (\et, \th, \r) \quad \hbox{ for every $ h \in L^2(Q)$,}
\end{align*}
where $(\et, \th, \r)$ is the unique solution
to system \EQLin\ corresponding to $h$ obtained from Theorem \ref{THM_Wp_linear}.
\Ethm

To derive an explicit expression from \eqref{optimal_formal} we are left
with the task of proving the \Frechet~differentiability of the reduced cost functional
$\Jred$ with respect to the time variable.
This can be performed rigorously by virtue of the improved regularity results obtained from Theorem~\ref{PROP_timereg}.
\Bthm
\label{THM_frechet_time}
Suppose that \accorpa{ab}{phzeropiccolo}~and \eqref{dtphzeroinV} hold, and 
in addition to \eqref{targets}, we assume \eqref{targets_reg}.
Moreover, let $(u,\tau)$ be an admissible control pair with the corresponding state $(\m,\ph,\s)$.
Then, the reduced
cost functional $\Jred$ is \Frechet~differentiable
with respect to time and
\begin{align}
	\non
	D_\tau\Jred({u},{\tau})
	&=
	\frac \bQ 2 \iO |\ph({\tau}) - \phQ({\tau})|^2
	+ \bO \iO (\ph({\tau}) - \phO) \dt \ph({\tau})
	\\ & \quad\quad
	+ \frac \bQh 2 \iO |\s({\tau}) - \sQ({\tau})|^2
	+ \frac \bmas 2 \iO \dt\ph({\tau})
	+ \btime
	+b_6 (\tau - \ts).
	\label{DJ_tau}
\end{align}
\Ethm
\Bdim
It readily follows from computing the derivative. Moreover, let us notice that
the terms
$ \bO \iO (\ph({\tau}) - \phO) \dt \ph({\tau})$ and $ \frac \bmas 2 \iO \dt\ph({\tau})$ are meaningful by virtue of the refined result Theorem \ref{THM_Wp_state_system}.
For more details we refer to \cite{GARLR}, where
the authors showed how the time derivative for time-dependent cost functionals can be 
performed in a general setting.
\Edim

\Bthm
\label{THM_Opt_cond_prima}
Assume that \accorpa{ab}{phzeropiccolo}~and \eqref{dtphzeroinV} are fulfilled.
Furthermore, in addition to \eqref{targets} we assume \eqref{targets_reg}, and let
$(\overline{u},\overline{\tau})$ be an optimal control for $\CP$ with the corresponding
optimal state $(\bm,\bph,\bs)$ obtained from Theorems \ref{THM_Wp_state_system}, and \ref{PROP_timereg}. Then, $(\overline{u},\overline{\tau})$
necessarily fulfils the following variational inequality
\begin{align}
	\non
	& \bQ \int_{Q_{\overline{\tau}}} (\bph-\phQ) \th 
	+ \bO \iO (\bph(\overline{\tau})-\phO) \th(\overline{\tau}) 
	+ \bQh  \int_{Q_{\overline{\tau}}} (\bs-\sQ) \r
	+ \frac \bmas 2 \iO \th(\overline{\tau}) 
	\\ & \quad\quad
	+ \bz \intQ \overline{u}(v-\overline{u})
	\geq
	0
	\quad \hbox{for every $v \in \Uad$},
	\label{foc_first}
\end{align}
where $(\eta,\th,\r)$ is the unique solutions to the linearised system \EQLin~corresponding to
$h=v-\overline{u}$ obtained from Theorem \ref{THM_Wp_linear}.
Moreover, we have that
\begin{align}
 D_\tau\Jred(\overline{u},\overline{\tau}) 
	\begin{cases}
	\geq 0 \quad \hbox{if $\overline{\tau}=0$,}
	\\
	= 0 \quad \hbox{if $\overline{\tau} \in (0,T)$},
	\\
	\leq 0 \quad \hbox{if $\overline{\tau}=T$,}
	\label{DJ_optim}
\end{cases}
\end{align}
where $D_\tau\Jred(\overline{u},\overline{\tau})$ 
is given by \eqref{DJ_tau} evaluated 
at the optimum pair $(\overline{u},\overline{\tau})$.
In addition, setting
\begin{align*}
	 \Lambda(\uopt, \tauopt)&:= 
	\frac \bQ 2 \iO |\ph({\tauopt}) - \phQ({\tauopt})|^2
	+ \bO \iO (\ph({\tauopt}) - \phO) \dt \ph({\tauopt})
	\\ & \quad
	+ \frac \bQh 2 \iO |\s({\tauopt}) - \sQ({\tauopt})|^2
	+ \frac \bmas 2 \iO \dt\ph({\tauopt})
	+ \btime,
\end{align*}
it follows that
$D_\tau\Jred(\overline{u},\overline{\tau}) = \Lambda(\uopt, \tauopt) 
+ b_6 (\tauopt - \ts)$. Hence, if $b_6\not = 0$,  
the condition \eqref{DJ_optim} can be implicitly characterised as
\begin{align}
	\begin{cases}
	\Lambda(\uopt, 0) \geq {b_6}\ts  \quad &\hbox{if $\overline{\tau}=0$,}
	\\
	\tauopt = \ts - {b_6}^{-1} \Lambda(\uopt, \tauopt) \quad &\hbox{if $\overline{\tau} \in (0,T)$,}
	\\
	\Lambda(\uopt, T) \leq  b_6 (\ts-T)  \quad &\hbox{if $\overline{\tau}=T$.}
	\label{DJ_optim_espl}
\end{cases}
\end{align}
\Ethm
\Bdim
As already mentioned, the variational inequalities \eqref{foc_first} and \eqref{DJ_optim} 
directly follows by exploiting the abstract conditions \eqref{optimal_formal}. 
As \eqref{foc_first} is concerned, 
let us notice that, loosely speaking, $\Jred$ is 
the composition of $\J$ with $\S$ so that it suffices to combine the \Frechet\
differentiability of the two operators with
the chain rule to get \eqref{foc_first}. In this direction, let us 
introduce the auxiliary operator
$\widetilde \S:\UR\to {\cal X} \times \UR$,
where 
\begin{align*}
	{\cal X}:=(\H1H \cap \L\infty V \cap \L2 W)^2,
\end{align*} 
defined by $\widetilde \S(u ):=(\S_2(u),\S_3(u),u)$
so that $\Jred(u,\tau)=\J \bigl(\widetilde\S(u),\tau \bigr)$.
Then, from Theorem~\ref{THM_Frechet_space} we infer that
\Beq
  D \widetilde\S(u):
  h \mapsto 
  (\th,\r,h)
  \quad \hbox{for every $h \in \UR$},
  \non
\Eeq
where $(\et,\th,\r)$ is the unique solution to the linearised system 
\EQLin\ corresponding to~$h$ obtained from Theorem \ref{THM_Wp_linear}.
Moreover, it is \sfw~to realise that $\J$ is \Frechet\ differentiable with respect to $u$
and that, for every $\tau\in [0,T]$, it holds that
\begin{align*}
  \non
  [D_u \J (\ph, \s, u, \tau)](\Phi, \Psi, h, \tau)
  & = \bQ \int_{Q_\tau} (\ph - \phQ) \Phi 
  + \bO \iO (\ph(\tau)-\phO) \Phi(\tau)
  \\
  & \quad  
  + \bQh \int_{Q_\tau} (\s - \sQ) \Psi 
  + \frac \bmas 2 \iO \Phi(\tau)
  + \bz \intQ u h
    \\
  & \quad 
  \hspace{2.3cm}
  \quad \hbox{for every $(\Phi, \Psi, h)\in {\cal X} \times \UR $}.
\end{align*}
Thus, we invoke the chain rule to obtain that
\begin{align*}
  [D_u\J_{\textrm red}(\uopt, \overline\tau)](h, \overline{\tau})
  &= [ D_u \J \bigl( \widetilde \S (\uopt),\overline\tau \bigr)] 
  	\bigl( [D\widetilde\S(\uopt)] (h), \overline\tau \bigr)
  \\ & 
  = [D_u \J(\bph,\bs,\uopt,\overline\tau)] (\th,\r,h,\overline\tau)
  \\ & \non
  =
  \bQ \int_{Q_{\overline{\tau}}} (\bph-\phQ) \th 
	+ \bO \iO (\bph(\overline{\tau})-\phO) \th(\overline{\tau}) 
	+ \bQh  \int_{Q_{\overline{\tau}}} (\bs-\sQ) \r
	\\ & \quad\quad
	+ \frac \bmas 2 \iO \th(\overline{\tau}) 
	+ \bz \intQ \overline{u}h,
\end{align*}
which leads to \eqref{foc_first}.

As for the second inequality \eqref{DJ_optim}, it readily follows from the second of
\eqref{optimal_formal} along with the characterisation given by \eqref{DJ_tau}.
Lastly, the first and the last conditions
of \eqref{DJ_optim} are consequences of the fact that we cannot 
exclude the cases $s=0$ and $s=T$, while the middle one
follows from the fact that, whenever $\overline{\tau}\in (0,T)$, 
we can simply take $s=\overline{\tau} \pm \zeta$,
with $\zeta>0$, to argue that $D_\tau\Jred(\overline{u},\overline{\tau}) = 0$.
\Edim

Let us emphasise that the above characterisation \eqref{DJ_optim_espl}
is new with respect to the previous contributions
\cite{CRW} and \cite{GARLR}, where just the condition \eqref{DJ_optim} was obtained.
In fact, the more explicit condition \eqref{DJ_optim_espl} can be now obtained by virtue of
the additional tracking-type term $\tfrac {b_6}2 |\tau-\ts|^2$ that we have added
in the cost functional.

\subsection{Adjoint System and First-order Optimality Condition}
\label{subsection}
This section is devoted to the introduction and discussion of the 
adjoint system to \EQ\
which is a key argument in simplifying the variational inequality \eqref{foc_first}.
Only \sfw~modifications are in order with respect to \cite{S} and
it can be easily shown that the (formal) {\it adjoint system} to \EQ\ reads as
\begin{align}
   -\b \dt q  - \dt p- \Delta q + F''(\bph)q + P'(\bph)(\bs - \bm)(r-p)&= \bQ(\bph - \phQ)
  \quad &&\hbox{in $ Q_{\overline{\tau}}$,}
  \label{EQAggprima}
  \\
   -\a \dt p - \Delta p - q + P(\bph)(p -r) &= 0
  \label{EQAggseconda}
  \quad &&\hbox{in $Q_{\overline{\tau}}$,}
  \\
   -\dt r - \Delta r + P(\bph)(r - p) &= \bQh (\bs - \sQ)
  \label{EQAggterza}
  \quad &&\hbox{in $Q_{\overline{\tau}}$,}
  \\
  \dn q = \dn p = \dn r & = 0
  \quad &&\hbox{on $\Sigma_{\overline{\tau}}$,}
  \label{BCEQAgg}
  \\
   \b q(\overline{\tau}) =\bO(\bph(\overline{\tau}) - \phO)+ {\tfrac \bmas2},  \ \ 
  p(\overline{\tau})& = 0, 
  \ \  r(\overline{\tau}) = 0 
  \quad &&\hbox{in $\Omega$.}
  \label{ICEQAgg}
\end{align}
\Accorpa\EQAgg EQAggprima ICEQAgg
The main differences with respect to \cite{S} are:
\begin{enumerate}
\item[(i)] The system \EQAgg\ has to be considered for $a.a. t \in (0,\overline{\tau})$
instead of $a.a. t \in (0,T)$.
\item[(ii)] In the final condition of $\b q(\overline{\tau})$ there appears a new term 
$ \tfrac \bmas 2$ which is due to the presence of $\frac \bmas 2 \iO (1+\ph(\tau))$
in the cost functional.
\item[(iii)] The final condition for $r$ is zero since
we are not considering any term involving the contribution 
$\iO |\s(\tau)-\sO|^2$ in the cost functional.
\end{enumerate}

It is worth noting that the above system is a backward-in-time problem with terminal data for
$q$ only belonging to $\Lx2$ (see assumption \eqref{targets}).
Therefore, the first equation \eqref{EQAggprima} has to be considered in a weak sense, i.e. we require that
\begin{align*}
	   & -{}_{V^*}\< \dt( p + \b q) (t),v >_V
	   +\iO \nabla q(t)\cdot \nabla v 
	   +\iO F''(\bph(t))q(t) v
 	\\ & \quad\qquad
 	+\iO P'(\bph(t))(\bs(t) - \bm(t))(r(t)-p(t))v
	 \\& \qquad
	 = \iO \bQ(\bph(t) - \phQ(t))v
	  \quad \hbox{for every $v\in V$ and, $a.a. t \in (0,\overline{\tau})$.}
\end{align*}

\Bthm
\label{THM_adjoint}
Assume that the assumptions \tutteleipotesi~are verified. 
Then, the adjoint system \EQAgg~admits a unique 
solution $(q,p,r)$ such that
\begin{align}
	q & \in H^1(0,\overline{\tau}; \Vp) \cap  L^{\infty}(0,\overline{\tau}; H)
	\cap L^2(0,\overline{\tau}; V)
	\subset C^0([0,\overline{\tau}]; H),
	\\ 
	p,r &\in H^1(0,\overline{\tau}; H) \cap L^{\infty}(0,\overline{\tau}; V)
	\cap L^{2}(0,\overline{\tau}; W) \subset C^0([0,\overline{\tau}]; V).
	\label{regadj}
\end{align}
\Ethm

\begin{proof}
As before, we proceed formally since the approach is standard and the system is linear.

\noindent
{\bf First estimate}
First, we add to both sides of \eqref{EQAggseconda} and \eqref{EQAggterza} 
the terms $p$ and $r$, respectively.
Then, we multiply \eqref{EQAggprima} by $q$, the new \eqref{EQAggseconda} by $\dt p$,
the new \eqref{EQAggterza} by $\dt r$, add the resulting equalities,
integrate over $Q_{\tauopt}$ 
and by parts. After some rearrangements
and a cancellation, we obtain that
\begin{align*}
	& \frac \b2 \iO|{q(t)}|^2
	+ \int_{Q_{\tauopt}} |\nabla q|^2
	+\a \int_{Q_{\tauopt}} |\dt p|^2
	+ \frac 12 \normaV{p(t)}^2
	+ \int_{Q_{\tauopt}} |\dt r|^2
	+ \frac 12 \normaV{r(t)}^2
	\\ & \quad 
	=
	 \frac \b2 \iO|{q(\tauopt)}|^2
	+ \int_{Q_{\tauopt}} \bQ(\bph - \phQ)q
	- \int_{Q_{\tauopt}}\bQh (\bs - \sQ)\dt r
	\\ & \qquad 
	- \int_{Q_{\tauopt}} F''(\bph) |q|^2
	- \int_{Q_{\tauopt}} P'(\bph)(\bs-\bm)(r-p)q
	+\int_{Q_{\tauopt}} P(\bph)(p-r)\dt p
	\\ & \qquad 
	- \int_{Q_{\tauopt}}p \, \dt p
	+\int_{Q_{\tauopt}} P(\bph)(r-p)\dt r
	-\int_{Q_{\tauopt}}r \, \dt r,
\end{align*}
where we denote by $I_1,...,I_{9}$ the integrals on the \rhs.
Using the Young inequality and recalling \eqref{initial_data} and \eqref{stimasep},
we easily obtain that
\begin{align*}
	&|I_1|+|I_2|+|I_3|+|I_4|+|I_7|+|I_{9}|
	\\ &\quad
	\leq 
	\delta \int_{Q_{\tauopt}} (|\dt p|^2+ |\dt r|^2)
	+ \cd \int_{Q_{\tauopt}} (|q|^2+|p|^2+|r|^2)
	+c,
\end{align*}
for a positive $\d$ yet to be determined and for 
a positive constant $\cd$ which only depends on $\d$.
Furthermore, using the \Holder\ and Young inequalities we infer that,
for every $\d>0$,
\begin{align*}
	|I_5|+|I_6|+|I_8|
	& \leq
	c  \int_{0}^{\tauopt} (\norma{\bs}_4+\norma{\bm}_4)(\norma{r}_4+\norma{p}_4)\norma{q}_2
	\\ & \quad 
	+  c \int_{Q_{\tauopt}} (|p|+|r|)(|\dt p|+|\dt r|)
	\\ & \leq 
	c  \int_{0}^{\tauopt} (\norma{\bs}_V+\norma{\bm}_V)(\norma{r}_V+\norma{p}_V)\norma{q}_H
	\\ & \quad 
	+  \d \int_{Q_{\tauopt}} (|\dt p|^2+|\dt r|^2)
	+ \cd \int_{Q_{\tauopt}} (|p|^2+|r|^2)
	\\ & \leq
	c  \int_{0}^{\tauopt} (\norma{\bs}_V^2+\norma{\bm}_V^2)(\norma{r}_V^2+\norma{p}_V^2)
	+ c \int_{Q_{\tauopt}} |q|^2
	\\ & \quad 
	+  \d \int_{Q_{\tauopt}} (|\dt p|^2+|\dt r|^2)
	+ \cd \int_{Q_{\tauopt}} (|p|^2+|r|^2),
\end{align*}
where we also used the continuous embedding $V \subset \Lx4$ and that $\bm$ and $\bs$,
as solutions to \EQ, satisfy \eqref{regstima} so that $(\norma{\bs}_V^2+\norma{\bm}_V^2) \in L^\infty (0,T)$.
Lastly, adjusting $\d\in(0,1)$ small enough, a Gronwall argument yields that
\begin{align*}
	\norma{q}_{L^{\infty}(0,\overline{\tau}; H)
	\cap L^2(0,\overline{\tau}; V)}
	+\norma{p}_{H^1(0,\overline{\tau}; H) \cap L^{\infty}(0,\overline{\tau}; V)}
	+\norma{r}_{H^1(0,\overline{\tau}; H) \cap L^{\infty}(0,\overline{\tau}; V)}
	\leq c.
\end{align*}

\noindent
{\bf Second estimate}
Next, comparison in \eqref{EQAggseconda} and then in 
\eqref{EQAggterza}, along with the above estimate, and 
the elliptic regularity theory allow us to deduce that
\begin{align*}
	\norma{p}_{L^{2}(0,\overline{\tau}; W)}
	+\norma{r}_{L^{2}(0,\overline{\tau}; W)}
	\leq c.
\end{align*}

\noindent
{\bf Third estimate}
Lastly, by comparison in \eqref{EQAggprima} we immediately realise that
\begin{align*}
	\norma{\dt q}_{L^{2}(0,\overline{\tau}; \Vp)}
	\leq c,
\end{align*}
which conclude the proof since the uniqueness directly follows from the above a priori estimates
by classical arguments since the adjoint system \EQAgg\ is linear in $(q,p,r)$.
\end{proof}

\an{
Here, let us point out that the terminal condition
of the adjoint variable $q$ given by \eqref{ICEQAgg}
allows us to rewrite the sum of the second and fourth terms of \eqref{DJ_tau}, i.e. $\bO \iO (\ph({\tau}) - \phO) \dt \ph({\tau})+\tfrac \bmas 2 \iO \dt\ph({\tau})$,
as $\b\iO \dt \ph(\tau)q(\tau)$.
This latter can be characterised also in different ways.
In fact, upon multiplying by $\b$ the weak formulation of \eqref{EQ_prima}
and the weak formulation of \eqref{EQ_seconda} by $q$, we infer that, for all $t\in[0,T]$,
\begin{align*}
	\b \iO  \dt \ph(t) q(t)
	= \iO \Delta \ph(t) q(t)
	- \iO F'(\ph(t)) q(t)
	+\iO \m(t)q(t),
\end{align*}
as well as
\begin{align*}
	\b \iO  \dt \ph(t) q(t)
	= \b \iO P(\ph(t))(\s(t)-\m(t)) q(t)
	- \a\b\iO \dt\m(t)q(t)
	+ \b \iO \Delta \m(t) q(t).
\end{align*}
These are completely meaningful in a pointwise sense if we can guarantee that
\begin{align*}
	 \s,q \in \C0 H,  \ \dt\m, \dt\ph \in \C0 {H}, \ \hbox{ and } \ \m,\ph \in \C0 {\Hx2}.
\end{align*}
Indeed these requirements are fulfilled under the framework of the refined regularity results expressed by Theorem \ref{PROP_timereg} and Theorem \ref{PROP_timeregdue}.}

Next, by using the adjoint variables $(q,p,r)$
we can eliminate the linearised variables from the variational 
inequality \eqref{foc_first} producing a simpler formulation for
the first-order necessary conditions for optimality.
\Bthm
\label{THM_final_FOC}
Suppose that \accorpa{ab}{phzeropiccolo}~and \eqref{dtphzeroinV} are satisfied.
Let $(\overline{u},\overline{\tau})$,
$(\bm,\bph,\bs)$ and $(q,p,r)$ 
be an optimal control pair with the corresponding
state and adjoint variables obtained from Theorems \ref{THM_Wp_state_system}, \ref{PROP_timereg}, and \ref{THM_adjoint}, respectively.
Then, $(\overline{u},\overline{\tau})$ necessarily satisfies
\Beq
	\label{final_foc_prima}
	\int_{Q_{\overline{\tau}}} r 	
	(v- \overline{u})
	+
	\bz \intQ  \overline{u}(v- \overline{u})
	\geq 0 \quad \hbox{for every $v\in\Uad$.}
\Eeq
\Ethm

\Bcor
Let the assumptions of Theorem \ref{THM_final_FOC} be satisfied.
Moreover, let us set $\widetilde{r}$ as the zero extension
of $r$ in $[0,T]$. 
Then, the optimal pair $(\overline{u},\overline{\tau})$ necessarily satisfies
\Beq
	\label{final_foc_simply}
	\intQ (\widetilde{r} +\bz \overline{u})(v - \overline{u}) \geq 0 \quad \hbox{for every $v\in\Uad$.}
\Eeq
Moreover, if $\bz\not =0$, the optimal control $\overline{u}$
is nothing but the 
$\L2 H$-orthogonal projection of $-{\bz}\!^{-1}\,{\widetilde{r}}$ onto the closed subspace $\Uad$.
\Ecor

\Brem
Note that the case $\overline{\tau}=0$ covers a 
special and trivia role. 
In fact, in this case 
the above variational inequality \eqref{final_foc_prima} reduces to
\begin{align}
	\non
	 \bz \intQ\overline{u}(v-\overline{u})
	\geq
	0
	\quad \hbox{for every $v \in \Uad$},
\end{align}
which, whenever $\bz>0$, yields that $\overline{u}$ is the orthogonal projection
of $0$ onto the closed subspace $\Uad$.
The same consequence can be drawn from evaluating the variational inequality \eqref{foc_first}
at $\overline{\tau}=0$ and using that $\th(0)=0$.
\Erem

\Brem
As a consequence of \eqref{final_foc_simply},
we can identify, via Riesz's representation
theorem, the gradient of the reduced cost functional
as
\Beq
	\non
	\nabla \Jred(\overline{u},\overline{\tau})=\widetilde{r} + \bz \overline{u}.
\Eeq
This fact is extremely important 
from the numerical viewpoint since it implies the possibility to analyse 
the optimal control problem $\CP$ 
as a constrained minimisation problem
via standard techniques (e.g. by applying the conjugate gradient method).
\Erem 

\begin{proof}[Proof of Theorem~\ref{THM_final_FOC}]
By virtue of simplicity we proceed formally by employing the 
formal Lagrangian method (see, e.g., \cite{Trol,Lions_OPT}).
Comparing the two variational inequalities \eqref{foc_first} and \eqref{final_foc_prima},
we realise that it suffices to check that
\begin{align}
	\non
	\int_{Q_{\overline{\tau}}} r h
	&=
	\bQ \int_{Q_{\overline{\tau}}} (\bph-\phQ) \th 
	+ \bO \iO (\bph(\overline{\tau})-\phO) \th(\overline{\tau}) 
	\\ &\quad 
	+ \bQh  \int_{Q_{\overline{\tau}}} (\bs-\sQ) \r
	+ \frac \bmas 2 \iO \th(\overline{\tau}) ,
	\label{condition}
\end{align}
where $h$ is taken as $h=v-\uopt$ and $(\et,\th,\r)$ is the unique solution to \EQLin\ associated to $h$ obtained from Theorem \ref{THM_Wp_linear}.
In this direction, we multiply \eqref{EQLinprima} by $p$, \eqref{EQLinseconda} by $q$,
\eqref{EQLinterza} by $r$, and integrate over $Q_{\overline{\tau}}$ to get
\begin{align*}
 	0  = &\int_{Q_{\overline{\tau}}} p 
 				[\a \dt \et + \dt \th - \Delta \et 
 				- P'(\bph) (\bs - \bm)\th - P(\bph)(\r - \et)]
 		  \\
 		  &+\int_{Q_{\overline{\tau}}} q
 		  		[\b \dt \th - \Delta \th + F''(\bph)\th-\et]
 		  \\
 		  &+\int_{Q_{\overline{\tau}}} r
 		  		[\dt \r - \Delta \r + P'(\bph) (\bs - \bm)\th + P(\bph)(\r - \et) - h].
\end{align*}
Then, we move the last term to the \lhs~and integrate by parts to obtain that
\begin{align*}
	\int_{Q_{\overline{\tau}}} r h
	&=
	\int_{Q_{\overline{\tau}}} \et 
 				[ -\a \dt p - \Delta p - q + P(\bph)(p -r) ]
 		  \\
 		  &\quad +\int_{Q_{\overline{\tau}}} \th
 		  		[-\b \dt q  - \dt p- \Delta q + F''(\bph)q + P'(\bph)(\bs - \bm)(r-p)]
 		  \\
 		  &\quad +\int_{Q_{\overline{\tau}}} \r
 		  		[-\dt r - \Delta r + P(\bph)(r - p)]
 		  \\
 		  &\quad +\iO [\a \et (\tauopt)p(\tauopt)+ \th(\tauopt)p(\tauopt)
 		  +\b \th(\tauopt)q(\tauopt)+\r(\tauopt)r(\tauopt)],
\end{align*}
where we also owe to the homogeneous Neumann boundary conditions for the
linearised and adjoint variables, and to the initial conditions for the
linearised variables. Finally, accounting for the adjoint system
\EQAgg, we conclude that the above equation reduces to
\eqref{condition}.
\end{proof}
\section{Some Possible Generalisations}
\label{SEC_GENERALIZING_COST_FZ}
\setcounter{equation}{0}
In the remainder of the paper, we aim at providing some 
indications concerning some further generalisations.
First, we will show how to possibly overcome the issue already
mentioned regarding the control of the nutrient variable $\s$ at the given time $\tau$.
Next, we will spend some words concerning a similar minimisation problem 
in which the role of the control variable sligthly differs from our choice.

\subsection{A Relaxation Argument}
From the mathematical viewpoint, a natural term to be considered in the cost functional
is $\iO |\s(\tau)-\sO|^2$.
However, as already emphasised, to give meaning to the necessary
conditions that will eventually appear, further temporal regularity for $\s$
has to be established. This will force us to demand 
the control $u$ to be more regular in time, say $\H1 H$,
which is unrealistic in the application.
Anyhow, a possible way to overcome this issue could be to follow
the relaxation strategy employed in \cite{GARLR}. 
To this aim, let us fix a positive constant $\eps$
and define the {\it relaxed} cost functional $\J_\eps$ as follows
\begin{align}
	 \non
	\J_\eps (\ph, \s, u, \tau)  &: = 
	\J(\ph, \s, u, \tau)
	+ \frac \bOh {2\eps} \int_{\tau-\eps}^{\tau} \iO |\s-\sO|^2,
    \label{costfunct_relaxed}
\end{align}
for a non-negative constant $\bOh$ and for a given target function $\sO$.
Note that the factor $\tfrac 1\eps$ is due to normalisation since
$\frac 1 {\eps} \int_{\tau-\eps}^{\tau} = 1$.
With this adjustement on the cost functional, we can control the final configuration
of the nutrient without demanding any additional regularity 
for the nutrient variable $\s$.
Hence, instad of considering $\CP$ we consider
\begin{align}
	\non
	&
	{(CP)_\eps} 
	\quad 
	\inf_{(\ph,\s,u, \tau) \in {\cal A}_{ad}}
	\J_\eps(\ph,\s,u, \tau).
\end{align}
The most part of the results follow in the same way.
Hence, we proceed schematically just mentioning the main differences.

\noindent
{\bf Existence}
The first arrangement to be done concerns the existence of a minimiser (cf. Section \ref{SEC_EX_MIN}).
The proof can be reproduced using the direct method of calculus of variations provided to
explain how the new term of the cost functional can be handled.
In this direction, let us point out that, along with \eqref{chi_tau_n}, 
we also have, as $n\to \infty$,
\Beq
	\non
	\chi_{[\tau_n-\eps,\tau_n]} (\cdot)
	\to 
	\chi_{[\overline{\tau}-\eps,\overline{\tau}]}(\cdot)
	\quad \hbox{$\aat$}.
\Eeq
Hence, by similar reasoning, we also conclude that, as $n\to \infty$,
\Beq
	\non
	\frac \bOh {2\eps} \int_{\tau_n-\eps}^{\tau_n} 
	\iO |\s_n({\tau}_n)-\sO|^2
	\to
	\frac \bOh {2\eps} \int_{\overline{\tau}-\eps}^{\overline{\tau}} 
	\iO |\bs({\overline{\tau}})-\sO|^2,
\Eeq
while the rest of the proof is exactly the same as in the proof of Theorem \ref{THM_Existence_optimal_controls}.

\noindent
{\bf \Frechet~differentiability of the reduced cost functional}
As expected, the main differences are related to the \Frechet~differentiability of
the corresponding reduced cost functional.
In fact, the corresponding of \eqref{foc_first} becomes
\begin{align*}
	\non
	& \bQ \int_0^{\overline{\tau}} \iO (\bph-\phQ) \th 
	+ \bO \iO (\bph(\overline{\tau})-\phO) \th(\overline{\tau}) 
	+ \bQh  \int_0^{\overline{\tau}} \iO (\bs-\sQ) \r
	\\ & \quad\quad
	+ \frac {\bOh} {\eps} \int_{\overline{\tau}-\eps}^{\overline{\tau}} \iO (\bs-\sO) \r
	+ \frac \bmas 2 \iO \th(\overline{\tau}) 
	+ \bz \intQ \overline{u}(v-\overline{u})
	\geq
	0
	\quad \hbox{for every $v \in \Uad$}.
\end{align*}
As the time derivative is concerned, we have to adjust a little the framework by 
assuming that $\sO \in H^1(-{\eps}, T; H) $ and that the variable $\s$ 
is meaningful for negative time. Hence, we simply
postulate that $\s(t) := \s_0$ if $t<0$.
Thus, the corresponding \Frechet\ derivative with respect to time reads as
\begin{align*}
	\non
	D_\tau\Jred({u},{\tau})
	&=
	\frac \bQ 2 \iO |\ph({\tau}) - \phQ({\tau})|^2
	+ \bO \iO (\ph({\tau}) - \phO) \dt \ph({\tau})
	\\ & \quad\quad \non
	+ \frac \bQh 2 \iO |\s({\tau}) - \sQ({\tau})|^2
	\\ & \quad\quad \non
	+ \frac {\bOh}{2\eps} \biggl( \iO |(\s- \sO)(\tau)|^2 - \iO |(\s- \sO)(\tau-{\eps})|^2 \biggr)
	\\ & \quad\quad  
	+ \frac \bmas 2 \iO \dt\ph({\tau})
	+ \btime
	+b_6 (\tau-\ts).
\end{align*}

\noindent
{\bf The adjoint system}
Lastly, the adjoint system slightly differs and becomes
\begin{align*}
  & -\b \dt q - \dt p - \Delta q + F''(\bph)q + P'(\bph)(\bs - \bm)(r-p)= \bQ(\bph - \phQ)
  \quad &&\hbox{in $\, Q_{\overline{\tau}},$}
  \\
  & -\a \dt p - \Delta p - q + P(\bph)(p -r) = 0
  \quad &&\hbox{in $\,Q_{\overline{\tau}},$}
  \\
  & -\dt r - \Delta r + P(\bph)(r - p) = \bQh (\bs - \sQ)
  +\tfrac \bOh{{\eps}} \chi_{(\overline{\tau}-{\eps},\overline{\tau})}(\cdot)\,  (\bs-\sO)  
  \quad &&\hbox{in $\,Q_{\overline{\tau}},$}
  \\
  &\dn q = \dn p = \dn r = 0
  \quad &&\hbox{on $\Sigma_{\overline{\tau}},$}
  \\
  & p(\overline{\tau}) + \b q(\overline{\tau}) =\bO(\bph(\overline{\tau}) - \phO)+ {\tfrac \bmas2},  \ \ 
  \a p(\overline{\tau}) = 0, 
  \ \  r(\overline{\tau}) = 0 
  \quad &&\hbox{in $\Omega$.}
\end{align*}
Notice that the only difference is the
\rhs~of the third equation which, however, still belongs to $\L2 H$
and therefore the well-posedness of the above system easily follows adapting the lines of argument of Theorem~\ref{THM_adjoint}.

\subsection{A Different Control Variable}
Let us conclude the paper by pointing out that another popular choice 
for the control variable $u$ is the one employed in \cite{GARLR,EK,EK_ADV}.
There, the control variable $u$ is placed in 
equation \eqref{EQ_prima} and it models the elimination of tumour cells by the effect of a 
cytotoxic drug. Namely, we can consider the following state system
\begin{align*}
   \a \dt \m + \dt \ph - \Delta \m &= P(\phi) (\sigma - \mu) - \kappa u  h(\ph) 
  \quad &&\hbox{in $\, Q$,}
  \\
   \mu &= \beta \dt \phi - \Delta \phi + F'(\phi)
  \quad &&\hbox{in $\,Q$,}
  \\
  \dt \sigma - \Delta \sigma  &= - P(\phi) (\sigma - \mu)
  \quad &&\hbox{in $\,Q$,}
  \\
   \dn \m &= \dn \ph =\dn \s = 0
  \quad &&\hbox{on $\,\Sigma$,}
  \\
   \m(0)&=\m_0,\, \ph(0)=\ph_0,\, \s(0)=\s_0
  \quad &&\hbox{in $\,\Omega,$}
\end{align*}
where $\kappa$ is a positive constant, whereas
the symbol $h$ stands for an interpolation function which
vanishes at $-1$ and attains value $1$ at $1$.
Moreover, the control $u$ ranges between $0$ and $1$ in order
to model no dosage and full dosage of the drug, respectively.
So, when $\ph=-1$
no drug is dispensed, when $\ph=1$ there is a full dosage of
the drug, and in between there is an intermediate supply.

It is worth noting that $\kappa h(\ph) u \in L^\infty(Q)$ so that the same arguments
employed in \cite{S} can be reproduced in the same manner to 
obtain the corresponding of Theorem \ref{THM_Wp_state_system}.
However, notice that Theorem \ref{PROP_timereg} cannot directly be reproduced.
In fact, since the control variable is now placed in the phase equation, 
it will be necessary to require that $u\in\H1 H$.
Therefore, in the spirit of the above section 
one is reduced to consider a relaxed cost functional
\begin{align*}
	 \non
	\J_\eps (\ph, \s, u, \tau)  &: = 
	\frac {a_1} 2 \int_{Q_\tau} |\ph-\phQ|^2
	+ \frac {a_2} {2\eps} \int_{\tau-\eps}^{\tau} \iO |\ph-\phO|^2
	+ \frac {a_3} 2 \int_{Q_\tau} |\s-\sQ|^2
	\\ \non
	& \quad\quad
	{+ \frac {a_4} {2}\iO |\s-\sO|^2}
	+ \frac {a_5}{2\eps} \int_{\tau-\eps}^{\tau} \iO (1+\ph)
	+ {a_6} \tau
	+ \frac {a_7}2 |\tau-\ts|^2
	+ \frac {a_0} {2} \intQ |u|^2,
    \label{}
\end{align*}
for some non-negative constants $a_0,...,a_7$. 
Let us claim that, providing to require some natural assumptions, 
the variable $\s$ may enjoy higher temporal regularity so that the third term
in the cost functional above can be considered without any 
relaxation arguments. 
Let us claim that, after proving higher temporal regularity for the variable $\s$,
the expected optimality conditions read as
\begin{align*}
	{a_7} \intQ  \uopt (v-\uopt)- \kappa \int_{Q_{\overline \tau}} h(\bph)p (v-\uopt)
	\geq 0 \quad
	\hbox{for every $v\in\Uad$,}
\end{align*}
where $p$ stands for the associated adjoint variable
and
\begin{align*}
	& 
	\frac {a_1}2 \iO |\ph(\overline \tau)-\phQ(\overline \tau)|^2
	+
	\frac {a_2}{2\eps} \iO \bigl(|(\ph - \phO)(\overline \tau)|^2-|(\ph - \phO)(\overline \tau - \eps)|^2\bigr)
	\\ & \quad
	+\frac {a_3}{2}\iO |\s(\overline \tau)-\sQ(\overline \tau)|^2
	+  {a_4} \iO (\s(\overline \tau)-\sO) \dt\bs(\overline \tau)
	\\ & \quad
	+\frac {a_5}{2\eps} \iO (\ph (\overline \tau)-\ph (\overline \tau - \eps))
	+ a_6
	+a_7 (\tau-\ts)
	\begin{cases}
	\geq 0 \quad \hbox{if $\overline{\tau}=0$,}
	\\
	= 0 \quad \hbox{if $\overline{\tau} \in (0,T)$,}
	\\
	\leq 0 \quad \hbox{if $\overline{\tau}=T$.}
	\end{cases}
\end{align*}
The details are left to the reader.

\small
\subsection*{Acknowledgments}
The author wishes to express his gratitude to Professor Pierluigi Colli for several useful 
discussions and suggestions which have improved the manuscript.
Moreover, the author would like to thank Professor Elisabetta 
Rocca whose useful comments have helped to clarify some technical points.
\vspace{3truemm}
\footnotesize

\Begin{thebibliography}{10}

\bibitem{Agosti}
A. Agosti, P.F. Antonietti, P. Ciarletta, M. Grasselli and M. Verani,
A Cahn--Hilliard--type equation with application to tumor growth dynamics, 
{\it Math. Methods Appl. Sci.}, {\bf 40} (2017), 7598-–7626.
%


\bibitem{BRZ}
H. Brezis,
``Op\'erateurs maximaux monotones et semi-groupes de contractions dans les
espaces de Hilbert'', North-Holland Math. Stud. {\bf 5}, North-Holland, Amsterdam, 1973.

\bibitem{CRW}
C. Cavaterra, E. Rocca and H. Wu,
Long--time Dynamics and Optimal Control of a Diffuse Interface Model for Tumor Growth,
{\it Appl. Math. Optim.} (2019), https://doi.org/10.1007/s00245-019-09562-5.


\bibitem{CGH}
P. Colli, G. Gilardi and D. Hilhorst,
On a Cahn--Hilliard type phase field system related to tumor growth,
{\it Discrete Contin. Dyn. Syst.} {\bf 35} (2015), 2423--2442.

\bibitem{SM}
P. Colli, G. Gilardi, G. Marinoschi and E. Rocca,
Sliding mode control for a phase field system related to tumor growth,
{\it Appl. Math. Optim.} to appear (2019), doi:10.1007/s00245-017-9451-z.
%

\bibitem{CGRS_VAN}
P. Colli, G. Gilardi, E. Rocca and J. Sprekels,
Vanishing viscosities and error estimate for a Cahn-–Hilliard type phase field system related to tumor growth,
{\it Nonlinear Anal. Real World Appl.} {\bf 26} (2015), 93--108.

\bibitem{CGRS_OPT}
P. Colli, G. Gilardi, E. Rocca and J. Sprekels,
Optimal distributed control of a diffuse interface model of tumor growth,
{\it Nonlinearity} {\bf 30} (2017), 2518--2546.

\bibitem{CGRS_ASY}
P. Colli, G. Gilardi, E. Rocca and J. Sprekels,
Asymptotic analyses and error estimates for a Cahn--Hilliard~type phase field system modeling tumor growth,
{\it Discrete Contin. Dyn. Syst. Ser. S} {\bf 10} (2017), 37--54.

\bibitem{CLLW}
V. Cristini, X. Li, J.S. Lowengrub, S.M. Wise,
Nonlinear simulations of solid tumor growth using a mixture model: invasion and branching.
{\it J. Math. Biol.} {\bf 58} (2009), 723–-763.

\bibitem{CL}
V. Cristini, J. Lowengrub,
Multiscale Modeling of Cancer: An Integrated Experimental and Mathematical
{\it Modeling Approach. Cambridge University Press}, Leiden (2010).

\bibitem{DFRGM}
M. Dai, E. Feireisl, E. Rocca, G. Schimperna, M. Schonbek,
Analysis of a diffuse interface model of multispecies tumor growth,
{\it Nonlinearity\/} {\bf  30} (2017), 1639.

\bibitem{EK_ADV}
M. Ebenbeck and P. Knopf,
Optimal control theory and advanced optimality conditions for a diffuse interface model of tumor growth,
{\it ESAIM: COCV} (2020) https://doi.org/10.1051/cocv/2019059.
 
\bibitem{EK}
M. Ebenbeck and P. Knopf,
Optimal medication for tumors modeled by a Cahn--Hilliard--Brinkman equation,
{\it Calc. Var.} {\bf 58} (2019) https://doi.org/10.1007/s00526-019-1579-z.

\bibitem{EGAR}
M. Ebenbeck and H. Garcke,
Analysis of a Cahn–-Hilliard–-Brinkman model for tumour growth with chemotaxis.
{\it J. Differential Equations,} (2018) https://doi.org/10.1016/j.jde.2018.10.045.

\bibitem{FGR}
S. Frigeri, M. Grasselli, E. Rocca,
On a diffuse interface model of tumor growth,
{\it  European J. Appl. Math.\/} {\bf 26 } (2015), 215--243. 

\bibitem{FLRS}
S. Frigeri, K.F. Lam, E. Rocca, G. Schimperna,
On a multi-species Cahn--Hilliard-Darcy tumor growth model with singular potentials,
{\it Comm. in Math. Sci.} {\bf  (16)(3)} (2018), 821--856. 

\bibitem{FRL}
S. Frigeri, K.F. Lam and E. Rocca,
On a diffuse interface model for tumour growth with non-local interactions and degenerate 
mobilities,
In {\sl  Solvability, Regularity, and Optimal Control of Boundary Value Problems for PDEs},
P. Colli, A. Favini, E. Rocca, G. Schimperna, J. Sprekels (ed.),
{\it Springer INdAM Series,} {\bf 22}, Springer, Cham, 2017.

%
\bibitem{GARL_3}
H. Garcke and K. F. Lam,
Global weak solutions and asymptotic limits of a Cahn--Hilliard--Darcy system modelling tumour growth,
{\it AIMS Mathematics} {\bf 1 (3)} (2016), 318--360.
\bibitem{GARL_1}
H. Garcke and K. F. Lam,
Well-posedness of a Cahn-–Hilliard–-Darcy system modelling tumour
growth with chemotaxis and active transport,
{\it European. J. Appl. Math.} {\bf 28 (2)} (2017), 284--316.
\bibitem{GARL_2}
H. Garcke and K. F. Lam,
Analysis of a Cahn--Hilliard system with non-zero Dirichlet 
conditions modeling tumor growth with chemotaxis,
{\it Discrete Contin. Dyn. Syst.} {\bf 37 (8)} (2017), 4277--4308.
\bibitem{GARL_4}
H. Garcke and K. F. Lam,
On a Cahn--Hilliard--Darcy system for tumour growth with solution dependent source terms, 
in {\sl Trends on Applications of Mathematics to Mechanics}, 
E.~Rocca, U.~Stefanelli, L.~Truskinovski, A.~Visintin~(ed.), 
{\it Springer INdAM Series} {\bf 27}, Springer, Cham, 2018, 243--264.
\bibitem{GAR}
H. Garcke, K. F. Lam, R. N\"urnberg and E. Sitka,
A multiphase Cahn--Hilliard--Darcy model for tumour growth with necrosis,
{\it Mathematical Models and Methods in Applied Sciences} {\bf 28 (3)} (2018), 525--577.

\bibitem{GARLR}
H. Garcke, K. F. Lam and E. Rocca,
Optimal control of treatment time in a diffuse interface model of tumor growth,
{\it Appl. Math. Optim.} {\bf 78}(3) (2018), {495--544}.

\bibitem{GLSS}
H. Garcke, K.F. Lam, E. Sitka, V. Styles,
A Cahn–-Hilliard-–Darcy model for tumour growth with chemotaxis and active transport.
{\it Math. Models Methods Appl. Sci. } {\bf 26(6)} (2016), 1095-–1148.


\bibitem{OHP}
A. Hawkins, J.T Oden, S. Prudhomme,
General diffuse-interface theories and an approach to predictive tumor growth modeling. 
{\it Math. Models Methods Appl. Sci.} {\bf 58} (2010), 723–-763. 

\bibitem{HDPZO}
A. Hawkins-Daarud, S. Prudhomme, K.G. van der Zee, J.T. Oden,
Bayesian calibration, validation, and uncertainty quantification of diffuse 
interface models of tumor growth. 
{\it J. Math. Biol.} {\bf 67} (2013), 1457–-1485. 

\bibitem{HDZO}
A. Hawkins-Daruud, K. G. van der Zee and J. T. Oden, Numerical simulation of
a thermodynamically consistent four-species tumor growth model, Int. J. Numer.
{\it Math. Biomed. Engng.} {\bf 28} (2011), 3–-24.

\bibitem{HKNZ}
D. Hilhorst, J. Kampmann, T. N. Nguyen and K. G. van der Zee, Formal asymptotic
limit of a diffuse-interface tumor-growth model, 
{\it Math. Models Methods Appl. Sci.} {\bf 25} (2015), 1011-–1043.



\bibitem{Lions_OPT}
J.-L. Lions,
Contr\^ole optimal de syst\`emes gouverne\'s par des equations aux d\'eriv\'ees partielles,
Dunod, Paris, 1968.



\bibitem{Mir_CH}
A. Miranville,
The Cahn--Hilliard equation and some of its variants, 
{\it AIMS Mathematics,} {\bf 2} (2017), 479–-544.

\bibitem{MRS}
A. Miranville, E. Rocca, and G. Schimperna,
On the long time behavior of a tumor growth model,
{\it Commun. Pure Appl. Anal.\/} {\bf 8} (2009) 881-912.
%

\bibitem{MirZel}
A. Miranville, S. Zelik, Attractors for dissipative partial differential equations in bounded and
unbounded domains, in “Handbook of Differential Equations: Evolutionary Equations, Vol. IV”
(eds. C.M. Dafermos and M. Pokorny), Elsevier/North-Holland, 103–200, 2008.

\bibitem{S_a}
A. Signori,
Vanishing parameter for an optimal control problem modeling tumor growth.
{\it Asymptot. Anal.} (2019), https://doi: 10.3233/ASY-191546.

\bibitem{S_b}
A. Signori,
Optimal treatment for a phase field system of Cahn--Hilliard 
type modeling tumor growth by asymptotic scheme, 
{\it Math. Control Relat. Fields} (2019), https://doi: 10.3934/mcrf.2019040.

\bibitem{S_DQ}
A. Signori,
Optimality conditions for an extended tumor growth model with 
double obstacle potential via deep quench approach, 
{\it Evol. Equ. Control Theory} (2019), https://doi: 10.3934/eect.2020003.

\bibitem{S}
A. Signori,
Optimal distributed control of an extended model of tumor
growth with logarithmic potential.
{\it Appl. Math. Optim.} (2018), https://doi.org/10.1007/s00245-018-9538-1.

\bibitem{Simon}
J. Simon,
{Compact sets in the space $L^p(0,T; B)$},
{\it Ann. Mat. Pura Appl.\/} 
{\bf 146~(4)} (1987) 65--96.

\bibitem{SW}
J. Sprekels and H. Wu,
Optimal Distributed Control of a Cahn-–Hilliard–-Darcy System with Mass Sources,
{\it Appl. Math. Optim.} (2019), https://doi.org/10.1007/s00245-019-09555-4.

\bibitem{Trol}
F. Tr\"oltzsch,
Optimal Control of Partial Differential Equations. Theory, Methods and Applications,
{\it Grad. Stud. in Math.,} Vol. {\bf 112}, AMS, Providence, RI, 2010.

\bibitem{WLFC}
S.M. Wise, J.S. Lowengrub, H.B. Frieboes, V. Cristini,
Three-dimensional multispecies nonlinear tumor growth—I: model and numerical method. 
{\it J. Theor. Biol.} {\bf 253(3)} (2008), 524–-543.

\bibitem{WZZ}
X. Wu, G.J. van Zwieten and K.G. van der Zee, Stabilized second-order splitting
schemes for Cahn--Hilliard~models with applications to diffuse-interface tumor-growth models, 
{\it Int. J. Numer. Meth. Biomed. Engng.} {\bf 30} (2014), 180--203.

\End{thebibliography}

\End{document}

\bye